\documentclass[review,onefignum,onetabnum]{siamart171218}

\usepackage[normalem]{ulem}
\usepackage{xcolor}
\usepackage{enumerate}
\usepackage{subcaption}

\usepackage{lipsum}
\usepackage{amsfonts}
\usepackage{graphicx}
\usepackage{epstopdf}
\usepackage{algorithmic}
\usepackage{comment}
\ifpdf
  \DeclareGraphicsExtensions{.eps,.pdf,.png,.jpg}
\else
  \DeclareGraphicsExtensions{.eps}
\fi

\newsiamremark{remark}{Remark}
\newsiamremark{hypothesis}{Hypothesis}
\crefname{hypothesis}{Hypothesis}{Hypotheses}
\newsiamthm{claim}{Claim}

\headers{Adjoint Least Squares Shadowing: Existence, Uniqueness and Coarse Domain Discretization}{P. Thakur and S. Nadarajah}

\title{Adjoint of Least Squares Shadowing: Existence, Uniqueness and Coarse Domain Discretization \thanks{Submitted to the editors.
\funding{This work was funded by the Fonds de recherche du Québec – Nature et technologies (FRQNT) Doctoral Award, McGill Engineering Doctoral Award (MEDA) and the Natural Sciences and Engineering Research Council of Canada (NSERC) Discovery Grant.}}}

\author{Pranshul Thakur \thanks{Department of Mechanical Engineering, McGill University, Montreal, H3A 0C3, QC, Canada (\email{pranshul.thakur@mail.mcgill.ca}).}
\and Siva Nadarajah \thanks{Department of Mechanical Engineering, McGill University, Montreal, H3A 0C3, QC, Canada (\email{siva.nadarajah@mcgill.ca}).} }

\usepackage{amsopn}
\usepackage{cancel}

\DeclareMathOperator{\diag}{diag}

\makeatletter
\newcommand*{\addFileDependency}[1]{
  \typeout{(#1)}
  \@addtofilelist{#1}
  \IfFileExists{#1}{}{\typeout{No file #1.}}
}
\makeatother

\newcommand{\vect}[1]{\boldsymbol{\mathbf{#1}}}

\ifpdf
\hypersetup{
  pdftitle={Adjoint of Least Squares Shadowing: Existence, Uniqueness and Coarse Domain Discretization},
  pdfauthor={P. Thakur, S. Nadarajah}
}
\fi

\begin{document}

\maketitle

\begin{abstract}
 Chaotic dynamical systems are characterized by the sensitive dependence of trajectories on initial conditions. Conventional sensitivity analysis of time-averaged functionals yields unbounded sensitivities when the simulation is chaotic. The least squares shadowing (LSS) is a popular approach to computing bounded sensitivities in the presence of chaotic dynamical systems. The current paper proves the existence, uniqueness, and boundedness of the adjoint of the LSS equations. In particular, the analysis yields a sharper bound on the condition number of the LSS equations than currently demonstrated in existing literature and shows that the condition number is bounded for large integration times. The derived bound on condition number also shows a relation between the conditioning of the LSS and the time dilation factor which is consistent with the trend numerically observed in the previous LSS literature. Furthermore, using the boundedness of the condition number for large integration times, we provide an alternate proof to~\cite{chater_wang_2017} of the convergence of the LSS sensitivity to the true sensitivity at the rate of $\mathcal{O}\left(\frac{1}{\sqrt{T}}\right)$ regardless of the boundary conditions imposed on the adjoint, as long as the adjoint boundary conditions are bounded. Existence and uniqueness of the solution to the continuous-in-time adjoint LSS equation ensure that the LSS equation can be discretized independently of the primal equation and that the true LSS adjoint solution is recovered as the time step is refined. This allows for the adjoint LSS equation to be discretized on a coarser time domain than that of the primal governing equation to reduce the cost of solving the linear space-time system.  
\end{abstract}



\section{Introduction}
\label{sec:introduction}
Sensitivity analysis of quantities of interest is used in many fields, such as optimization \cite{nocedal_wright}, computing improved functional estimates \cite{becker_rannacher_2001, giles_suli, pierce_giles}, mesh adaptation \cite{corrigan_et_al_2019,zahr1,thakur_nadarajah_2024, fidkowski_darmofal_2011, inriagoaloriented, dolejsigoal1} etc. The adjoint method of computing sensitivities is efficient when there are a lot fewer functionals than independent parameters, as is common in aerodynamic applications \cite{jameson_adjoint_optimization}. However, both the conventional forward and adjoint methods yield unbounded sensitivities when applied to chaotic systems \cite{lea_e_al_2000}. The reason is due to the positive Lyapunov exponents in chaotic systems, which cause most nearby trajectories to diverge exponentially. This phenomenon is commonly known as the `butterfly effect'. 

Several methods have been developed to compute the sensitivities of chaotic dynamical systems. The work of Ruelle \cite{ruelle_diff_of_srb} established the differentiability of chaotic hyperbolic systems, along with an equation to compute the linear response of the time-averaged functionals. Other approaches use the probability density function obtained from the Fokker-Plank equation \cite{thuburn_fokker_planck, blonigan_fokker_planck}. However, these methods rely on discretizing the entire phase space, which might be expensive when the number of degrees of freedom of the system is large. While it is possible to reduce the cost of probability distribution-based methods by discretizing only the region of the attractor instead of the entire phase space, such an approach requires the knowledge of the location of the attractor beforehand \cite{blonigan_fokker_planck}, which might be infeasible for cases involving high degrees of freedom. On the other hand, there are notable approaches that do not rely on probability distribution. The ensemble adjoint approach of Lea et al. \cite{lea_ensemble_adjoint, lea_e_al_2000} computes conventional adjoint solutions over short trajectories and takes the ensemble average of the resulting short-window sensitivities. However, the approach is known to converge slower than the Monte-Carlo approach to the correct sensitivity \cite{lea_ensemble_adjoint, wang_ensemble_adjoint}. It is known that hyperbolic systems with a compact attractor have trajectories that shadow the reference trajectory when the function generating the dynamics is perturbed \cite{anosov_shadowing, bowen_shadowing,numerical_shadowing,chater_wang_2017}. Using this idea, the work of Wang et al. \cite{wang_2013,wang_lss} developed a method of linearizing the governing equation using the shadowing trajectory. The shadow trajectory was computed either through splitting the governing equation into the invariant stable, unstable, and neutral subspaces \cite{wang_2013} or by solving a least-squares minimization problem, yielding the least-squares shadowing approach (LSS) \cite{wang_lss}. In LSS, a space-time linear system needs to be solved, which requires significant memory and computational resources for large test cases \cite{blonigan_airfoil}. To reduce the cost of the LSS approach, there have been approaches that reduce the size of the linear system through multiple shooting \cite{blonigan_mss}, reduced-order modeling \cite{vermier_lss, fidkowski_lss} or formulating the shadowing algorithm in a Fourier space with truncation \cite{ashley_hicken_fourier,kantarakias_papadakis}. Some approaches replace the initial condition with periodic boundary conditions on the tangent/adjoint equations to yield bounded sensitivities, such as the unstable periodic orbits \cite{lasagna_upo} and periodic shadowing \cite{lasagna_periodic_shadowing}. Recently, the unstable periodic orbit approach was extended to discontinuous periodic orbits to yield regularized but accurate responses to parameter perturbations \cite{hicken_chaos}. The non-intrusive \cite{nilss} and adjoint non-intrusive methods \cite{blonigan_adjoint_nilss,nilsas} have been developed which require only the unstable subspace of the dynamical system for forward sensitivity and an additional one-dimensional neutral subspace for adjoint-based sensitivity \cite{nilsas} methods. Their cost scales with the number of unstable Lyapunov covariant vectors in the system and are efficient for problems with lower dimensional unstable subspace, but are expensive when the number of positive Lyapunov exponents is large. We also note that the shadowing-based approaches might compute shadowing trajectories that are non-physical \cite{chandramoorthy_wang_nonphysical} and approaches based on the linear response \cite{ruelle_diff_of_srb} have been investigated \cite{chandramoorthy_efficient_linear_response,ni2023recursive,sliwiak2023} to address the issue. The current paper does not aim to address this issue as it is not clear if the problem occurs for cases involving a large number of degrees of freedom.

In \cite{wang_gomez_lss}, a condition number for the LSS equation was derived and it was noted to be similar to the condition number of the Poisson equation. The current work shows that the condition number of LSS is, in fact, bounded as the integration time $T$ increases. This fact is then used to establish the existence and uniqueness of the LSS solution and leads to an alternate proof of the convergence of the LSS sensitivity. The current paper is organized as follows: \Cref{sec:notations} describes the notations used in the paper. \Cref{sec:lss_introduction} provides a brief introduction to the least-squares shadowing (LSS) method. \Cref{sec:adjoint_lss_introduction} provides the adjoint of LSS in the form used later in the proofs. \Cref{sec:adjoint_shadowing_direction} gives a brief overview of the adjoint shadowing direction, which is used as the true reference adjoint solution. \Cref{sec:convergence_adjoint_lss} proves that the LSS condition number is bounded for a long integration time. This result is then used to provide an alternative proof of the convergence of LSS sensitivity to the true sensitivity with an adjoint formulation. \Cref{sec:convergence_adjoint_lss} also shows that the Sobolev norm of the error in adjoint is bounded, the boundary conditions of the adjoint LSS do not affect the convergence of the scheme and the error in sensitivity due to discretization is of the order of the local truncation error. \cref{rm:alpha_squared_dependence} then compares the effect of the time dilation factor of LSS on the derived bounds for condition number and the trend is found to be consistent with that numerically observed in previous literature on LSS \cite{blonigan_multigrid}. Using the result of the existence and uniqueness of the solution to the adjoint LSS, the continuous adjoint LSS is discretized on a coarser time domain than that used to compute the primal solution trajectory, resulting in a smaller space-time linear system. \Cref{sec:numerical_cases} presents numerical test cases of applying the adjoint LSS for sensitivity analysis.      
 
\section{Preliminaries}
\label{sec:notations}
Consider the following dynamical system, governing the evolution of a state with time $t$:
\begin{equation}
\label{eq:governing_eq}
    \frac{du}{dt} = f(u(t,s),s),\;\;u(0) = u_0
\end{equation}
where $u(t,s) \in X$ is the state in a $n$-dimensional Hilbert space $X$ with an initial condition $u_0$, $s \in \mathbb{R}$ is the design or control parameter (for example, the shape of an airfoil \cite{jameson_adjoint_optimization} or uncertainty parameters to determine robust flight trajectory of a re-entry vehicle \cite{optimal_control}) and $f(u,s)$ is a smooth function governing the system. The solution trajectory is used to compute an output, $J(u,s) : X \times \mathbb{R} \rightarrow \mathbb{R}$. The quantity of interest is a time-averaged output: 

\begin{equation}
\label{eq:time_averaged_functional}
    \bar{J}(s) = \lim_{T\to \infty} \frac{1}{T} \int_0^T J(u,s) dt.
\end{equation}
We assume that the system is ergodic \cite{katok_dynamical_system}. Therefore, $\bar{J}$ does not depend on the initial condition $u_0$ and is only a function of $s$.

Perturbing $s$ modifies the nonlinear function $f$ governing the system. Thus, the resulting perturbed trajectory $\tilde{u}$ is different from $u$ and, for chaotic dynamical systems, $\tilde{u}$ diverges exponentially from $u$ in the vicinity of $u$ for most initial conditions $\tilde{u}(0)$. This is due to the positive Lyapunov exponent of the chaotic systems. Conventional sensitivity analysis assumes that the perturbed system has the same initial condition as the unperturbed system, i.e. $\tilde{u}(0) = u(0)$, leading to unbounded growth of the norm of $v = du/ds$ after linearization. On the other hand, the true sensitivity
\begin{equation}
    \label{eq:true_sensitivity}
    \frac{d\bar{J}}{ds} = \lim_{\delta s \to 0} \frac{\bar{J}(s+\delta s) - \bar{J}(s)}{\delta s}
\end{equation}
exists and is bounded for hyperbolic chaotic systems \cite{ruelle_diff_of_srb}.

According to the shadowing lemma \cite[Theorem 18.1.2]{katok_dynamical_system}, \cite{bowen_shadowing,anosov_shadowing,pilyugin_shadowing}, there is an initial condition $\tilde{u}(0)$ and a time transformation $\tau(t,s)$ such that the trajectory of the perturbed system,
\begin{equation}
    \label{eq:perturbed_eq}
    \frac{d\tilde{u}}{d\tau} = f(\tilde{u}(\tau), s+\delta s),
\end{equation}
shadows the reference trajectory $u(t)$ $\forall t$; where $\tilde{u}$ is the solution to the perturbed governing equation.

Throughout the paper, we assume that the dynamical system \cref{eq:governing_eq} satisfies the assumption of uniform hyperbolicity \cite{katok_dynamical_system}. This means that the attractor can be split into invariant stable, unstable, and neutral subspaces. Further, there exists a set of Lyapunov exponents, $\{ \lambda_i\}_{i=1}^n$, governing the rate at which an infinitesimal perturbation in the state grows or shrinks with time. For each Lyapunov exponent $\lambda_i$, there is an associated covariant Lyapunov vector (CLV), $\phi_i(u)$, satisfying \cite{ginelli_2007,ruelle_ergodic,eckmann_ergodic, pilyugin_shadowing,wang_2013}:
\begin{equation}
\label{eq:forward_clv}
    \frac{d\phi_i}{dt} = f_u\phi_i - \lambda_i \phi_i.
\end{equation}
Thus, a perturbation in the state along $\phi_i$ grows exponentially at the rate $e^{\lambda_i t}$. The vectors $\{\phi_i(u)\}_{i=1}^n$ are bounded at all times and form a basis of phase space at each point on the attractor. Similarly, whenever CLVs exist, Oseledet's multiplicative ergodic theorem \cite{oseledets} also implies the existence of adjoint covariant Lyapunov vectors, $\{\hat{\phi}_i\}_{i=1}^n$, as shown in \cite{adjoint_clv}. The adjoint covariant vectors satisfy \cite{wang_2013, ni_adjoint_arxiv,ginelli_2007}:
\begin{equation}
\label{eq:adjoint_clv}
    -\frac{d\hat{\phi}_i}{dt} = f_u^*\hat{\phi}_i - \lambda_i \hat{\phi}_i.
\end{equation}
Like covariant Lyapunov vectors, the adjoint covariant Lyapunov vectors are always bounded and form a basis of the phase space. Further, by choosing a proper normalization, the Lyapunov covariant and adjoint covariant vectors satisfy $\phi_i(t)^T\hat{\phi_j}(t) = \delta_{ij}$ \cite{ni_adjoint_arxiv,wang_2013}. By the assumption of uniform hyperbolicity, the neutral subspace is one-dimensional. Denoting $i=n_0$ to represent the neutral subspace, $\lambda_{n_0} = 0$, it can be seen from \cref{eq:forward_clv} that the neutral CLV is $\phi_{n_0} = f$ \cite{wang_2013,blonigan_multigrid}.

\section{Least squares shadowing approach to sensitivity analysis}
\label{sec:lss_introduction}
 In order to compute the true sensitivity, the least squares shadowing (LSS) approach of Wang et al. \cite{wang_lss} linearizes the trajectory with respect to the shadowing trajectory. Taking $\delta u(t) = \tilde{u}(\tau) - u(t)$ and linearizing the temporal transformation as $\frac{d\tau}{dt} = 1 + \eta(t)\delta s$, subtracting \cref{eq:governing_eq} from \cref{eq:perturbed_eq} yields \cite{wang_lss}
\begin{equation*}
    \frac{d \delta u}{dt} = f_u \delta u + f_s \delta s + \eta f \delta s
\end{equation*}
or 
\begin{equation}
      \label{eq:linearized_eq}
    \frac{d v}{dt} = f_u v + f_s + \eta f,
\end{equation}
where $v = \frac{du}{ds}$. In order to find the shadowing direction, the following least-squares minimization problem is solved:
\begin{equation}
    \label{eq:LSS_minimization}
    \begin{aligned}
        & \min_{v,\eta} & & \frac{1}{2}\int_0^T v^Tv + \alpha^2 \eta^2 dt \\
        & \textrm{s.t.} & & \frac{d v}{dt} = f_u v + f_s + \eta f,
    \end{aligned}
\end{equation}
with $\alpha^2 \in \mathbb{R}^+$. Note that \cref{eq:LSS_minimization} does not specify an initial condition. Instead, it solves an optimization problem to find a trajectory. The sensitivity of the functional \cref{eq:time_averaged_functional} can then be computed as:
\begin{equation}
\label{eq:functional_sensitivity}
    \frac{d \bar{J}}{ds} = \frac{1}{T} \int_0^T J_u v + J_s + \eta(J-\bar{J}) dt.
\end{equation}
Details of the derivation of this equation can be found in \cite{wang_lss,blonigan_multigrid}. It was shown in \cite{chater_wang_2017} that the sensitivity from \cref{eq:functional_sensitivity} with the shadowing direction from \cref{eq:LSS_minimization} converges to the true sensitivity in the limit as $T \to \infty$.

\section{Adjoint of least squares shadowing}
\label{sec:adjoint_lss_introduction} 
Solution to the convex minimization problem \cref{eq:LSS_minimization} can be obtained by formulating the Lagrangian $\mathcal{\tilde{L}}(v,\eta,w) = \frac{1}{2}\int_0^T v^Tv + \alpha^2 \eta^2 dt - \int_0^T w^T\left(\frac{d v}{dt} - f_u v - f_s - \eta f \right) dt$ and identifying $(v,\eta,w)$ such that the variation of $\mathcal{\tilde{L}}$ with respect to arbitrary perturbations $\delta v$, $\delta \eta$ and $\delta w$ vanishes. The resulting equation is a boundary value problem in time with the following system of equations \cite{wang_lss}:
\begin{equation}
\label{eq:forward_lss}
    \begin{aligned}
        & \frac{dw}{dt} + f_u^* w + v = 0,  \quad \quad w(0) = 0,\;\;w(T) = 0,\\
        & w^Tf = -\alpha^2\eta, \\
        & \frac{dv}{dt} - f_u v -f_s - \eta f = 0.
    \end{aligned}
\end{equation}

The adjoint of least squares shadowing can be derived by taking into account the sensitivity of the functional \cref{eq:functional_sensitivity}. Using Lagrange multipliers $r$, $\psi_2$ and $\psi$, the sensitivity can be computed as:
\begin{equation}
\label{eq:lagrangemult}
\begin{aligned}
    \frac{d \bar{J}}{ds} =  & \frac{1}{T} \int_0^T J_u v + J_s + \eta (J - \bar{J}) dt\\
    &- \frac{1}{T} \int_0^T r^T \left( \frac{dw}{dt} + f_u^* w + v  \right)dt\\
    & - \frac{1}{T} \int_0^T \psi_2 \left( w^Tf + \alpha^2\eta   \right)dt\\
    & - \frac{1}{T} \int_0^T \psi^T \left( \frac{dv}{dt} - f_u v -f_s - \eta f   \right)dt, 
\end{aligned}
\end{equation}
where the last three terms on the right-hand-side are zero from \cref{eq:forward_lss}. Performing integration-by-parts to eliminate $v$, $\eta$ and $w$ yields the sensitivity:
\begin{equation}
    \label{eq:adjoint_sensitivity}
    \frac{d \bar{J}}{ds} = \lim_{T\to\infty} \frac{1}{T} \int_0^T \left(J_s + \psi^Tf_s \right) dt,
\end{equation}
with the adjoint equations
\begin{equation}
\label{eq:adjoint_eq_pre}
\begin{aligned}
    & \frac{d \psi}{dt} + f_u^* \psi + J_u = r, \quad \quad \psi(0) = 0,\;\;\psi(T) = 0, \\
    & \alpha^2 \psi_2 = \psi^T f + J-\bar{J}, \\
    & \frac{dr}{dt} - f_u r = \psi_2 f.
\end{aligned}
\end{equation}
We note that the discrete version of the continuous adjoint LSS \cref{eq:adjoint_eq_pre} was derived in \cite{wang_lss}. By eliminating $\psi_2$, \cref{eq:adjoint_eq_pre} can be expressed as
\begin{subequations}
\label{eq:lss_adjoint}
\begin{align}
    & \frac{d \psi}{dt} + f_u^* \psi + J_u = r, \quad \quad \psi(0) = 0,\;\;\psi(T) = 0 \label{eq:lss_adjoint_a}\\
    & \frac{dr}{dt} - f_u r = \frac{1}{\alpha^2} \left(\psi^T f + J - \bar{J} \right) f.  \label{eq:lss_adjoint_b}
\end{align}
\end{subequations}

The next section presents the adjoint shadowing direction, which is taken as the true reference adjoint solution and is used to determine bounds on the solution to the adjoint LSS \cref{eq:lss_adjoint}.  

\section{Adjoint shadowing direction}
\label{sec:adjoint_shadowing_direction}
The conventional adjoint equation, given by
\begin{equation}
\label{eq:conventional_adjoint}
    \frac{d\psi^\infty}{dt} + f_u^* \psi^\infty + J_u = 0,
\end{equation}
has a solution that is unbounded for most terminal conditions. However, for uniformly hyperbolic chaotic systems, there is an adjoint solution, $\psi^\infty$, which satisfies \cref{eq:conventional_adjoint}, is bounded for all time $t \in [0,\infty)$ and satisfies \cite{wang_2013,ni_adjoint_arxiv} 
\begin{equation}
\label{eq:avg_adjoint}
    \frac{1}{T}\int_0^T \psi^{\infty T}f = 0,
\end{equation}
where $\psi^\infty$ is referred to as the adjoint shadowing direction.

It has been shown in \cite{wang_2013,ni_adjoint_arxiv} that the true sensitivity of the functional can be computed using the adjoint shadowing direction as
\begin{equation}
\label{eq:djds_shadowing}
    \frac{d \bar{J}^\infty}{ds} = \lim_{T\to\infty} \frac{1}{T} \int_0^T \left(J_s + \psi^{\infty T}f_s \right) dt,
\end{equation}
where the superscript $\infty$ in $\frac{d \bar{J}^\infty}{ds}$ refers to the use of $\psi^\infty$ in evaluating the sensitivity. For a detailed discussion of the adjoint shadowing direction, the reader is referred to \cite{ni_adjoint_arxiv}. Next, we present a result that would be useful to bound the adjoint LSS error (cf. \cite{ni_adjoint_arxiv, wang_2013}):
\begin{lemma}
\label{lm:adjoint_dot_f}
The shadowing direction, $\psi^\infty$, satisfying \cref{eq:conventional_adjoint} and \cref{eq:avg_adjoint} also satisfies
\begin{displaymath}
        \psi^{\infty}(t)^Tf(t) = \bar{J} - J(t) \;\; \forall t.
\end{displaymath}
\end{lemma}
\begin{proof}
    Taking the inner product of \cref{eq:conventional_adjoint} with $f(u,s)$ yields
    \begin{equation}
    \label{eq:psi_f_mid}
        f^T\left(\frac{d\psi^\infty}{dt} + f_u^* \psi^\infty + J_u \right) = \frac{d (\psi^{\infty T} f + J)}{dt} = 0 \implies \psi^{\infty T} f + J = C,
    \end{equation}
    where $C$ is a constant and we have used $\frac{df}{dt} = f_u \frac{du}{dt} = f_u f$ and $\frac{dJ}{dt} = J_u f$. Taking the time average on either side of the equation and using \cref{eq:avg_adjoint} yields
    \begin{equation*}
        \frac{1}{T} \int_0^T \psi^{\infty T}f dt + \bar{J} = C \implies C = \bar{J}
    \end{equation*}
    One obtains the stated equation after substituting $C$ in \cref{eq:psi_f_mid}.
\end{proof}

\section{Convergence of the adjoint LSS sensitivity}
\label{sec:convergence_adjoint_lss}
Convergence of the forward LSS sensitivity to the true sensitivity was shown in \cite{chater_wang_2017}. In this section, we provide an alternate proof of convergence of LSS sensitivity using the adjoint-based formulation, which leads to some additional insights on the adjoint solution. To that end, we take $\psi^{\infty}$ to be the true adjoint solution and denote the error in the adjoint as $e(t) = \psi(t) - \psi^{\infty}(t)$. 

\begin{lemma}
    \label{lm:error_bvp}
    The error $e(t)$ satisfies
    \begin{equation}
    \label{eq:operator_Le}
    \begin{aligned}
      \mathcal{L}e =  -\frac{d^2 e}{dt^2} & - \frac{d}{dt}(f_u^*e) + f_u \frac{de}{dt} + f_uf_u^*e + \frac{ff^T}{\alpha^2}e= 0 \\
     & e(0) = -\psi^{\infty}(0),\;\; e(T) = -\psi^{\infty}(T)
     \end{aligned}
\end{equation}
\end{lemma}
\begin{proof}
Subtracting \cref{eq:conventional_adjoint} from \cref{eq:lss_adjoint_a} and using \cref{lm:adjoint_dot_f} in \cref{eq:lss_adjoint_b} to note that 
\begin{equation*}
\psi^Tf + J - \bar{J} = \psi^Tf - \psi^{\infty T}f = e^Tf    
\end{equation*}
 yields
\begin{equation}
    \label{eq:error_bvp_1}
    \begin{aligned}
    & \frac{d e}{dt} + f_u^* e = r, \quad \quad e(0) = -\psi^{\infty}(0),\;\; e(T) = -\psi^{\infty}(T) \\
    & \frac{dr}{dt} - f_u r = \frac{1}{\alpha^2} \left(e^T f \right) f \\
\end{aligned}
\end{equation}

Eliminating $r(t)$ from the above equation yields the stated result.
\end{proof}
Note that $e(0)$ and $e(T)$ are bounded since $\psi^\infty$ is bounded at all times. Moreover, the LSS operator $\mathcal{L}$ in \cref{lm:error_bvp} is the same as that in \cite{wang_gomez_lss, blonigan_airfoil} and is known to be symmetric positive definite. It has been known that the absolute condition number of LSS, defined as $\|\mathcal{L}^{-1}\|$, is similar to that of Poisson's equation i.e it is proportional to $T^2$ \cite{wang_gomez_lss}. However, in the following section, we show that the condition number is in fact bounded for large integration times.

For the ease of derivations that follow, we express the error $e(t)$ as a non-homogeneous equation with homogeneous boundary conditions along with an added correction. Specifically, let
\begin{equation}
    e(t) = e_h(t) + e_p(t)
\end{equation}
satisfy the boundary conditions
    \begin{equation}
    \label{eq:ep_bc}
    \begin{aligned}
        & e_h(0) = e_h(T) = 0, \\
        & e_p(0) = e(0)\quad e_p(T) = e(T).
    \end{aligned}
    \end{equation}
The function $e_p(t)$ is chosen such that it satisfies the non-homogeneous boundary conditions specified in \cref{eq:ep_bc}. Using the operator $\mathcal{L}$ from \cref{lm:error_bvp},
\begin{equation}
\label{eq:L_eh_ep}
     \mathcal{L}e = 0 \implies \mathcal{L}e_h = - \mathcal{L}e_p
\end{equation}
Next, let
\begin{equation*}
\begin{aligned}
      & a = \left(\frac{e(0) - \exp(-T)e(T)}{1-\exp(-2T)} \right),\\
      & b = \left(\frac{e(T) - \exp(-T)e(0)}{1-\exp(-2T)} \right),
\end{aligned}
\end{equation*}
and
\begin{equation}
\label{eq:ep_exp}
    e_p(t) = a\exp(-t) + b\exp(-(T-t)).
\end{equation}
It can be verified that \cref{eq:ep_exp} satisfies the boundary conditions of $e_p(t)$ in \cref{eq:ep_bc}. 

Taking the inner product of \cref{lm:error_bvp} with $e_h(t)$ and using \cref{eq:L_eh_ep}, one obtains
\begin{equation}
\label{eq:bilinear_eh_ep_innerproduct}
    \mathcal{B}(e_h,e_h) = (e_h,\mathcal{L}e_h) = -(e_h,\mathcal{L}e_p).
\end{equation}
Using integration by parts with homogeneous boundary conditions of $e_h(t)$ yields (c.f. \cite{wang_gomez_lss})
\begin{equation}
\label{eq:eq_bilinear}
    \mathcal{B}(e_h,e_h) = \int_0^T \left( \frac{de_h}{dt} + f_u^* e_h \right)^T\left( \frac{de_h}{dt} + f_u^* e_h \right) dt  + \int_0^T \frac{1}{\alpha^2}e_h^T ff^Te_h dt.
\end{equation}

\begin{definition}
  We define the $L_2$, $H^1$ and the $W-$norms as
\begin{equation}
\label{eq:various_norms}
    \begin{aligned}
        & \|f\|_2 := \sqrt{\int_0^T f^Tf dt} \\
        & \|f\|_{H^1} := \sqrt{\|f\|^2_2 + \bigg\|\frac{df}{dt}\bigg\|^2_2}\\
        & \|f\|_W := \sqrt{\|f\|_2^2 + \bigg\|\frac{df}{dt} + f_u^*f \bigg\|_2^2 }
    \end{aligned}
\end{equation}

We also denote $\|f_u^*\|_{\infty} := \sup\limits_{t\in[0\;T]} \{\|f_u^*(t)\|_X\}$, where $\| f_u^*(t) \|_X = \sup\limits_{\substack{v\in X \\ \|v\|\leq 1}} \| f_u^*(t)v \|_X$ is the norm associated with the linear operator in the phase space $X$.  
\end{definition}

\begin{lemma}
\label{lm:norm_equivalence}
$\|f\|_{H^1} \leq (2 + \|f_u^*\|_{\infty})\|f\|_W$
\end{lemma}
\begin{proof}
    It can be seen that
    \begin{equation}
    \label{eq:ineq}
        \|f\|_W \geq \bigg\|\frac{df}{dt} + f_u^*f \bigg\|_2 \geq \bigg\|\frac{df}{dt}\bigg\|_2 - \|f_u^*\|_{\infty} \|f\|_2, 
    \end{equation}
        \begin{equation*}
      \implies \bigg\|\frac{df}{dt}\bigg\|_2 \leq  \|f\|_W +  \|f_u^*\|_{\infty} \|f\|_2 \leq (1 + \|f_u^*\|_{\infty})\|f\|_W,
    \end{equation*}
    where the second inequality in \cref{eq:ineq} is obtained from the reverse triangular inequality. Thus,
    \begin{equation*}
       \|f\|_2 + \bigg\|\frac{df}{dt}\bigg\|_2 \leq \|f\|_W + \bigg\|\frac{df}{dt}\bigg\|_2 \leq (2 + \|f_u^*\|_{\infty})\|f\|_W.
    \end{equation*}
    Using the fact that $\|f\|_{H^1} \leq  \|f\|_2 + \|\frac{df}{dt}\|_2$ yields the stated inequality.
\end{proof}

The next Lemma, \cref{lm:e_hat_bounds}, follows largely from Wang et al. \cite[Appendix A]{wang_gomez_lss} with the modification that the Poincar\'e-type inequality to bound the time derivative is derived without assuming the existence of Fourier series expansion. We include proofs of these results here for completeness.  

\begin{lemma}
\label{lm:e_hat_bounds}
Let $e_h(t) = \sum\limits_{i=1}^n \hat{e}_{hi}(t) \hat{\phi}_i(u(t))$, where $\{\lambda_i \}_{i=1}^n$ and $\{\hat{\phi}_i \}_{i=1}^n$ are the Lyapunov exponents and adjoint Lyapunov covariant vectors respectively. Then,

        \begin{equation} \tag{i}
        \label{eq:e_hat_inequality}
            c_{\min} \sum\limits_{i=1}^n \int_0^T \hat{e}_{hi}^2dt \leq \int_0^T e_h^Te_h dt \leq c_{\max} \sum\limits_{i=1}^n \int_0^T \hat{e}_{hi}^2 dt.
        \end{equation}
        \begin{equation} \tag{ii}
        \label{eq:dedt_bound_inequality}
         \sum_{i=1}^n c_{\min}\left(\frac{1}{T^2} + \lambda_i^2\right)\int_0^T \hat{e}_{hi}^2 dt \leq  \int_0^T \left( \frac{de_h}{dt} + f_u^* e_h \right)^T\left( \frac{de_h}{dt} + f_u^* e_h \right) dt.
         \end{equation}
    
with $0<c_{\min}<c_{\max}$.
\end{lemma}
\begin{proof}
    Since $\{\hat{\phi}_i(t) \}_{i=1}^n$ form a basis of the phase space \cite{ginelli_2007,ruelle_ergodic,eckmann_ergodic}, any $e_h(t)\in X$ can be expressed as $e_h(t) = \sum\limits_{i=1}^n \hat{e}_{hi}(t)\hat{\phi}_i(t)=\vect{\hat{\phi} \hat{e}_h}$. 
Adjoint covariant vectors, $\hat{\phi}_i (u(t))$, are bounded at all times \cite{ginelli_2007}. Therefore, $\vect{\hat{\phi}^T\hat{\phi}} \leq c_{\max}I$. The fact that $\vect{\hat{\phi}^T\hat{\phi}}$ is positive definite can be seen by noting that $\vect{\hat{\phi}^T\hat{\phi}}$ is symmetric, $\vect{v}^T \left(\vect{\hat{\phi}^T\hat{\phi}}\right)\vect{v} = \left(\vect{\hat{\phi}}\vect{v}\right)^T\left(\vect{\hat{\phi}}\vect{v}\right) = 0$ only for $\vect{v}=0$ because the columns of $\vect{\hat{\phi}}$ form a basis of phase space and  $\vect{v}^T \left(\vect{\hat{\phi}^T\hat{\phi}}\right)\vect{v}>0$ for $\vect{v}\neq \vect{0}$. Therefore $\exists\; c_{\min}>0$ such that $ c_{\min}I \leq \vect{\hat{\phi}^T\hat{\phi}}$. Noting that $e_h(t)^Te_h(t) = \vect{\hat{e}_h^T\hat{\phi}^T\hat{\phi}\hat{e}_h}$ and combining the previous two statements yields:
\begin{equation}
\label{eq:d_inequality}
  c_{\min} \sum\limits_{i=1}^n \hat{e}^2_{hi}(t)  \leq e_h(t)^Te_h(t) \leq c_{\max} \sum\limits_{i=1}^n \hat{e}^2_{hi}(t) \quad \forall t.
\end{equation}
Integrating \cref{eq:d_inequality} yields \cref{eq:e_hat_inequality}. 

Using the expansion $e_h(t) = \sum\limits_{i=1}^n \hat{e}_{hi}(t) \hat{\phi}_i(u(t))$ along with \cref{eq:adjoint_clv}, one obtains 
\begin{equation}
\label{eq:de_dt_hyperbolic}
     \frac{de_h}{dt} + f_u^* e_h  = \sum_{i=1}^n \left(\frac{d\hat{e}_{hi}}{dt} + \lambda_i \hat{e}_{hi} \right) \hat{\phi_i}.
\end{equation}

Substituting \cref{eq:de_dt_hyperbolic} on the right hand side of \cref{eq:dedt_bound_inequality} and using the bounds on $\vect{\hat{\phi}^T\hat{\phi}}$ yields:

\begin{equation}
     c_{\min}   \sum_{i=1}^n \int_0^T \left( \frac{d\hat{e}_{hi}}{dt} + \lambda_i \hat{e}_{hi} \right)^2 dt \leq \bigg\|\frac{de_h}{dt} + f_u^* e_h\bigg\|_2^2
\end{equation}

Expanding the summand on the left-hand-side, one obtains
\begin{equation}
\label{eq:expanded_dehat_dt}
     c_{\min}   \sum_{i=1}^n \int_0^T  \left(\frac{d\hat{e}_{hi}}{dt}\right)^2 dt + 2\lambda_i \cancel{\int_0^T \frac{d\hat{e}_{hi}}{dt}\hat{e}_{hi}dt} + \int_0^T\lambda_i^2 \hat{e}_{hi}^2 dt \leq \bigg\|\frac{de_h}{dt} + f_u^* e_h\bigg\|_2^2, 
\end{equation}
where the second term is zero because of the homogeneous boundary conditions of $e_h$.

Since  $\hat{e}_{hi}(0) = 0$, by the fundamental theorem of calculus, $\hat{e}_{hi}(t) = \int_0^t \frac{d\hat{e}_{hi}}{d\xi}d\xi$. Therefore, 
\begin{equation*}
    \hat{e}_{hi}(t) \leq \int_0^t \bigg| \frac{d\hat{e}_{hi}}{d\xi}\bigg| d\xi \leq \int_0^T \bigg| \frac{d\hat{e}_{hi}}{d\xi}\bigg|. 1 d\xi  \leq \sqrt{T} \bigg\| \frac{d\hat{e}_{hi}(\xi)}{d\xi}\bigg\|_2
\end{equation*}
\begin{equation}
\label{eq:poincare_inequality}
    \implies  \frac{1}{T^2} \int_0^T \hat{e}_{hi}^2 dt \leq \int_0^T \left( \frac{d\hat{e}_{hi}}{dt}\right)^2dt
\end{equation}
Substituting \cref{eq:poincare_inequality} into \cref{eq:expanded_dehat_dt} and noting that $\lambda_i$ is constant $\forall t$ yields \cref{eq:dedt_bound_inequality}.

\end{proof}

Next, we present a result that shows that the bilinear form is coercive as $T\to\infty$. 
\begin{lemma}
\label{lm:L_inverse_bound}
   $\exists\; 0<\gamma^2<\infty$ such that  $\gamma^2 \|e_h\|_{H^1}^2\leq (e_h,\mathcal{L}e_h)$ as $T \to \infty$. Specifically,
   \begin{equation*}
       \gamma^2 = \frac{1}{(2+\|f_u^*\|_\infty)^2}\min\bigg\{\frac{1}{2},  \frac{c_{\min} \{\lambda_i^2\}_{\substack{\\i=1,i\neq n_0}}^n}{2c_{\max}}, \frac{1}{\alpha^2 c_{\max}}   \bigg\}
   \end{equation*}
   is a constant independent of $T$. Moreover, $\|\mathcal{L}^{-1}\| \leq 1/\gamma^2$ is bounded.
\end{lemma}
\begin{proof}
    The bilinear form in \cref{eq:eq_bilinear} can be written as
    \begin{equation*}
        \mathcal{B}(e_h,e_h) = (e_h,\mathcal{L}e_h) = \mathcal{B}_1(e_h,e_h) + \mathcal{B}_2(e_h,e_h),
    \end{equation*}
    where
    \begin{equation*}
        \begin{aligned}
            & \mathcal{B}_1(e_h,e_h) = \frac{1}{2} \bigg\|\frac{de_h}{dt} + f_u^* e_h\bigg\|_2^2 \\
            & \mathcal{B}_2(e_h,e_h) = \frac{1}{2} \bigg\|\frac{de_h}{dt} + f_u^* e_h\bigg\|_2^2 + \int_0^T \frac{1}{\alpha^2}e_h^T ff^Te_h dt.
        \end{aligned}
    \end{equation*}

Representing $e_h(t)$ in terms of adjoint Lyapunov covariant vectors, we have \newline $e_h(t) = \sum\limits_{i=1}^n \hat{e}_{hi}(t)\hat{\phi}_i(t)$. Using \cref{eq:dedt_bound_inequality} in \cref{lm:e_hat_bounds} for the first term of $\mathcal{B}_2$ yields
\begin{equation}
  \sum_{i=1}^n \frac{c_{\min}}{2}\left(\frac{1}{T^2} + \lambda_i^2\right)\int_0^T \hat{e}_{hi}^2 dt + \int_0^T \frac{1}{\alpha^2}e_h^T ff^Te_h dt  \leq \mathcal{B}_2(e_h,e_h)
\end{equation}

Using the orthogonality property of adjoint covariant vectors, $e_h^Tf = \hat{e}_{n_0}(t)$. Therefore, we have
\begin{equation}
\label{eq:ef_to_e_hat}
     \sum_{i=1}^n \frac{c_{\min}}{2T^2} \int_0^T \hat{e}_{hi}^2 dt + \int_0^T \left(\sum_{i=1}^n \left( \frac{c_{\min}\lambda_i^2}{2} \hat{e}_{hi}^2\right) +  \frac{\hat{e}_{n_0}^2} {\alpha^2} \right)dt  \leq \mathcal{B}_2(e_h,e_h)
\end{equation}

Moreover, since $\lambda_{n_0}=0$, the integrand of the second term on the left-hand-side of \cref{eq:ef_to_e_hat} can be bounded from below by
\begin{equation}
    \label{eq:e_hat_take_min}
   \min \bigg\{ \frac{c_{\min} \{\lambda_i^2\}_{\substack{\\i=1,i\neq n_0}}^n}{2}, \frac{1}{\alpha^2} \bigg\} \sum_{i=1}^n \hat{e}_{hi}^2  \leq \sum_{i=1}^n \left( \frac{c_{\min}\lambda_i^2}{2} \hat{e}_{hi}^2\right) +  \frac{\hat{e}_{n_0}^2} {\alpha^2}
\end{equation}
The constant on the left-hand-side of \cref{eq:e_hat_take_min} is positive since, by the assumption of uniform hyperbolicity, $\lambda_{n_0}=0$ is the only zero Lyapunov exponent.

Substituting \cref{eq:e_hat_take_min} into \cref{eq:ef_to_e_hat} yields
\begin{equation*}
   \left(\frac{c_{\min}}{2T^2} + \min \bigg\{ \frac{c_{\min} \{\lambda_i^2\}_{\substack{\\i=1,i\neq n_0}}^n}{2}, \frac{1}{\alpha^2} \bigg\}\right) \int_0^T \sum_{i=1}^n \hat{e}_{hi}^2dt  \leq \mathcal{B}_2(e_h,e_h)
\end{equation*}
Using the inequality \cref{eq:e_hat_inequality} in \cref{lm:e_hat_bounds} yields
\begin{equation}
\label{eq:inequality_B2}
   \left(\frac{c_{\min}}{2c_{\max}T^2} + \min \bigg\{ \frac{c_{\min} \{\lambda_i^2\}_{\substack{\\i=1,i\neq n_0}}^n}{2c_{\max}}, \frac{1}{\alpha^2c_{\max}} \bigg\}\right) \|e_h\|_2^2 \leq \mathcal{B}_2(e_h,e_h)
\end{equation}
Since  $\mathcal{B}(e_h,e_h) = \mathcal{B}_1(e_h,e_h) + \mathcal{B}_2(e_h,e_h)$, adding $\mathcal{B}_1$ to \cref{eq:inequality_B2} and invoking the norm $\|\cdot\|_W$ in \cref{eq:various_norms} yields
\begin{equation*}
    \min \bigg\{ \frac{1}{2},  \left(\frac{c_{\min}}{2c_{\max}T^2} + \min \bigg\{ \frac{c_{\min} \{\lambda_i^2\}_{\substack{\\i=1,i\neq n_0}}^n}{2c_{\max}}, \frac{1}{\alpha^2c_{\max}} \bigg\}\right)\bigg\} \|e_h\|_W^2     \leq \mathcal{B}(e_h,e_h)
\end{equation*}
Using the norm equivalence relation \cref{lm:norm_equivalence},
\begin{equation*}
    \min \bigg\{ \frac{1}{2},  \frac{c_{\min}}{2c_{\max}T^2} + \min \bigg\{ \frac{c_{\min} \{\lambda_i^2\}_{\substack{\\i=1,i\neq n_0}}^n}{2c_{\max}}, \frac{1}{\alpha^2c_{\max}} \bigg\}\bigg\} \frac{\|e_h\|_{H^1}^2}{(2+\|f_u^*\|_\infty)^2}     \leq \mathcal{B}(e_h,e_h)
\end{equation*}
\begin{equation*}
    \xrightarrow{T\to\infty} \gamma^2 \|e_h\|^2_{H^1} \leq  \mathcal{B}(e_h,e_h)
\end{equation*}
as stated.

For an equation $\mathcal{L}f=g$, we have $f=\mathcal{L}^{-1}g$. Moreover,  $$\gamma^2\|f\|^2_{H^1} \leq (f,\mathcal{L}f)_2 = (f,g)_2 \leq \|f\|_2\|g\|_2 \leq  \|f\|_{H^1}\|g\|_{H^1}.$$ Hence, $\|f\|_{H^1} = \|\mathcal{L}^{-1}g\|_{H^1}  \leq \frac{\|g\|_{H^1}}{\gamma^2}$. Therefore, $\|\mathcal{L}^{-1}\| = \sup\limits_{\|g\|_{H^1}=1} \|\mathcal{L}^{-1}g\|_{H^1} \leq \frac{1}{\gamma^2}$ is bounded.
\end{proof}

\begin{remark}
\label{rm:L_inverse_l2}
    For an equation $\mathcal{L}f=g$, the derivation in \cref{lm:L_inverse_bound} also show that $\|f\|_2 \leq 1/\gamma^2 \|g\|_2$. Hence, the bound on the absolute condition number, $\|\mathcal{L}^{-1}\|$, is also valid in the $L_2$-norm.
\end{remark}

\begin{remark}
    It was shown in \cite{wang_gomez_lss} that $\|\mathcal{L}^{-1}\| < C T^2$, where the constant $C>0$. While the bound in \cite{wang_gomez_lss} indicates that $\|\mathcal{L}^{-1}\|$ grows quadratically with $T$ and can become unbounded, \cref{lm:L_inverse_bound} yields a sharper bound. In particular, it shows that $\|\mathcal{L}^{-1}\|$ is bounded as $T\to\infty$, allowing to establish the existence and uniqueness of the solution to the adjoint LSS.
\end{remark}

Consider the Sobolev space of functions:
\begin{equation}
    H^1 = \{f:[0\;T]\to X \;:\; \|f\|_{H^1} < \infty \}.
\end{equation}

\begin{corollary}
\label{cor:eh_is_in_sobolev}
    $e_h(t) \in H^1$.
\end{corollary}
\begin{proof}
        The right hand side of \cref{eq:L_eh_ep} can be expressed as
    \begin{equation*}
        -\mathcal{L}e_p = \frac{d^2 e_p}{dt^2} + \frac{d}{dt}(f_u^*e_p) - f_u \frac{de_p}{dt} - f_uf_u^*e_p - \frac{ff^T}{\alpha^2}e_p.
    \end{equation*}
    Substituting $e_p(t)$ from \cref{eq:ep_exp} yields
\begin{equation*}
    -\mathcal{L}e_p =  \exp(-t)h(t) + \exp(-(T-t))q(t)
\end{equation*}
where
\begin{equation*}
\begin{split}
        & h(t) = a + \frac{d}{dt}(f_u^*) a - f_u^*a + f_ua - f_uf_u^*a -\frac{1}{\alpha^2} ff^Ta \\
        & q(t) =  b + \frac{d}{dt}(f_u^*) b + f_u^*b - f_ub - f_uf_u^*b -\frac{1}{\alpha^2} ff^Tb \\
\end{split}
\end{equation*}

Assuming $f \in C^2(X)$ and the phase space $X$ is compact, $\|h\|_{\infty}$ and $\|q\|_{\infty}$ are bounded. Therefore
\begin{equation*}
\begin{split}
     \lim_{T\rightarrow \infty}\int_0^T (\mathcal{L}e_p)^T(\mathcal{L}e_p) dt &= \lim_{T\rightarrow \infty} \int_0^T \left( \exp(-2t)h^Th + 2\exp(-T)h^Tq + \exp(-2(T-t))q^Tq \right) dt \\ 
&\leq \lim_{T\rightarrow \infty} \left( \|h\|_{\infty}^2 \int_0^T \exp(-2t) dt + 2\sup\limits_{t\in[0 \;T]}\{h^Tq\} T\exp(-T) + \|q\|_{\infty}^2 \int_0^T \exp(-2(T-t)) dt \right) \\
& = \lim_{T\rightarrow \infty} \left(  \left(\|h\|_{\infty}^2 + \|q\|_{\infty}^2 \right) \frac{1 - \exp(-2T)}{2} + 2\max_t\{ h^Tq\} T\exp(-T) \right) \\
& =\frac{\|h\|_{\infty}^2 + \|q\|_{\infty}^2 }{2} 
\end{split}
\end{equation*}

which is a bounded constant. Moreover, using \cref{lm:L_inverse_bound} yields
\begin{equation*}
    \gamma^2\|e_h\|^2_{H^1} \leq (e_h,\mathcal{L}e_h)_2 = -(e_h,\mathcal{L}e_p)_2 \leq \|e_h\|_2\|\mathcal{L}e_p\|_2 \leq  \|e_h\|_{H^1}\|\mathcal{L}e_p\|_{2}.
\end{equation*}
\begin{equation*}
    \implies \lim_{T\to\infty} \|e_h\|_{H^1} \leq \frac{\|\mathcal{L}e_p\|_{2}}{\gamma^2} < \infty  
\end{equation*}
Hence, $e_h\in H^1$.
\end{proof}

The next result provides an upper bound on the LSS operator $\mathcal{L}$.
\begin{corollary}
\label{cor:bounded_bilinear}
    $\exists \;0<\beta^2 < \infty$ such that $$(\mathcal{L}w,v) \leq \beta^2\|w\|_{H^1}\|v\|_{H^1}\;\forall w,v \in H^1,$$ with $$ \beta^2 = (1 + \|f_u^*\|_{\infty})^2 + \frac{\|ff^T\|_{\infty}}{\alpha^2} .$$
\end{corollary}
\begin{proof}
For any $w,v \in H^1$, we have
\begin{align}
    (\mathcal{L}w,v) = & \int_0^T \left( \frac{dw}{dt} + f_u^* w \right)^T\left( \frac{dv}{dt} + f_u^* v \right) dt  + \int_0^T \frac{1}{\alpha^2}w^T ff^Tv dt \\
    & \leq \bigg\|\frac{dw}{dt} + f_u^* w \bigg\|_2 \bigg\|\frac{dv}{dt} + f_u^* v\bigg\|_2 + \frac{\|ff^T\|_\infty}{\alpha^2} \|w\|_2\|v\|_2 \\
    & \leq \left( (1 + \|f_u^*\|_{\infty})^2 + \frac{\|ff^T\|_{\infty}}{\alpha^2} \right) \|w\|_{H^1}\|v\|_{H^1}
\end{align}
where we have used the triangular inequality $\|\frac{dw}{dt}+f_u^* w \|_2 \leq \|\frac{dw}{dt} \|_2 + \| f_u^* w\|_2$ and fact that $\max\{\|w\|_2,\|\frac{dw}{dt}\|_2\} \leq \|w\|_{H^1}$.
\end{proof}

\begin{theorem}
    \label{thm:existence_uniqueness}
    The solution to the adjoint LSS equation \cref{eq:lss_adjoint} exists and is unique.
\end{theorem}
\begin{proof}
The solution $\psi$ to \cref{eq:lss_adjoint} can be expressed as $\psi = \psi^{\infty} + e_h(t) + e_p(t)$. $\psi^{\infty}(t)$ exists from the theory of dynamical systems \cite{ni_adjoint_arxiv,nilsas} while $e_p(t)$ exists by the choice of \cref{eq:ep_exp}. It was shown in \cref{cor:eh_is_in_sobolev} that, if $e_h$ exists, then $e_h\in H^1$. The Sobolev space $H^1$ is known to be a Hilbert space \cite{brenner_scott}, and is hence a complete metric space. Moreover, it was shown in \cref{lm:L_inverse_bound} that the bilinear form $\mathcal{B}(\cdot,\cdot)$ is coercive, while \cref{cor:bounded_bilinear} shows that $\mathcal{B}$ is bounded in $H^1$. Hence, the Lax-Milgram theorem \cite[Theorem 2.7.7]{brenner_scott} establishes the existence of a weak solution $e_h$ to the equation $\mathcal{L}e_h=-\mathcal{L}e_p$. 

Uniqueness follows from the fact that if two LSS adjoint solutions $\psi_1$ and $\psi_2$ satisfy \cref{eq:lss_adjoint}, then they also satisfy $\mathcal{L}(\psi_1-\psi_2)=0$. Therefore
\begin{equation*}
    \gamma^2 \|\psi_1-\psi_2\|^2_{H^1} \leq (\psi_1-\psi_2,\mathcal{L}(\psi_1-\psi_2)) = 0
\end{equation*}
\begin{equation*}
    \implies \|\psi_1-\psi_2\|_{H^1} = 0.
\end{equation*}
since $\gamma^2>0$. Thus, $\psi_1(t)=\psi_2(t)$ almost everywhere.  
\end{proof}

\begin{remark}
\label{rm:alpha_squared_dependence}
    It was observed numerically in Blonigan and Wang \cite{blonigan_multigrid} that the condition number of LSS depends on the time dilation factor $\alpha^2$ and deteriorates for large or small $\alpha^2$. The bounds obtained in \cref{lm:L_inverse_bound} and \cref{cor:bounded_bilinear} explain this trend. In particular, the relative condition number in the $H^1$-norm is
    \begin{equation*}
        \kappa_{\text{rel}} = \|\mathcal{L}^{-1}\|\|\mathcal{L}\| \leq \beta^2/\gamma^2.
    \end{equation*}
As $\alpha^2$ increases, $\gamma^2 \to \frac{1}{\alpha^2 c_{\max}(2+\|f_u^*\|_\infty)^2}$ in \cref{lm:L_inverse_bound} while $\beta^2 \to (1+\|f_u^*\|_\infty)^2$. Hence, $\kappa_{\text{rel}}$ increases and the conditioning deteriorates. On the other hand, by decreasing $\alpha^2$, $\beta^2$ in \cref{cor:bounded_bilinear} increases while $\gamma^2$ in \cref{lm:L_inverse_bound} is a constant, yielding large $\kappa_{\text{rel}}$. The values of $\gamma^2$ and $\beta^2$ indicate that the optimal $\alpha^2 \sim \mathcal{O}\left(\frac{1}{c_{\min} \min\{\lambda_i^2\}_{\substack{\\i=1,i \neq n_0}}^n}\right)$.  
\end{remark}

\begin{theorem}
\label{thm:sensitivity_convergence_rate}
The sensitivity $\frac{d\bar{J}}{ds}$ from the adjoint LSS converges to the true sensitivity at $\mathcal{O}\left(\frac{1}{\sqrt{T}}\right)$.
\end{theorem}
\begin{proof}
    Let $\frac{d \bar{J}}{ds}|_{\text{true}}$ be the true sensitivity, $\frac{d \bar{J}^\infty}{ds}$ the sensitivity computed using the adjoint shadowing direction $\psi^\infty$ \cref{eq:djds_shadowing} and $\frac{d \bar{J}}{ds}$ the sensitivity from adjoint LSS \cref{eq:adjoint_sensitivity}. Error in the LSS sensitivity can be expressed as 
\begin{equation}
\label{eq:sensitivity_error_total}
   \bigg| \frac{d\bar{J}}{ds} - \frac{d\bar{J}}{ds}\bigg|_{\text{true}} \bigg|  \leq \bigg|\frac{d \bar{J}}{ds} - \frac{d \bar{J}^\infty}{ds}\bigg| + \bigg|\frac{d \bar{J}^\infty}{ds} - \frac{d \bar{J}}{ds}\bigg|_{\text{true}} \bigg| 
\end{equation}
To bound the first term, we subtract \cref{eq:djds_shadowing} from \cref{eq:adjoint_sensitivity} to obtain
\begin{equation}
\label{eq:sensitivity_error_first_term}
  \bigg| \frac{d \bar{J}}{ds} - \frac{d \bar{J}^\infty}{ds} \bigg| < \bigg|\lim_{T\to\infty} \frac{1}{T}\int_0^T e^Tf_s dt \bigg| =  \bigg|\lim_{T\to\infty} \frac{1}{T}\int_0^T e_h^Tf_s dt + \lim_{T\to\infty} \frac{1}{T}\int_0^T e_p^Tf_s dt \bigg|
\end{equation}
The first term on the right-hand side can be expressed as
\begin{equation}
\label{eq:errbound1}
    \frac{1}{T}\int_0^T e_h^Tf_s dt \leq \frac{1}{T}\|e_h\|_2 \|f_s\|_2 \leq \frac{1}{T}\|e_h\|_2 \|f_s\|_{\infty} \sqrt{T} =  \frac{\|e_h\|_2 \|f_s\|_{\infty}}{\sqrt{T}}, 
\end{equation}
while the second term in \cref{eq:sensitivity_error_first_term} as
\begin{equation}
\label{eq:errbound2}
    \frac{1}{T}\int_0^T e_p^Tf_s dt \leq \frac{1}{T}\int_0^T \|e_p\|\|f_s\| dt \leq \frac{\|f_s\|_\infty \int_0^T \|e_p(t)\| dt}{T},
\end{equation}
with $\|e_p(t)\| = \sqrt{e_p(t)^Te_p(t)}$. Thus, using \cref{eq:errbound1,eq:errbound2} and the fact that for a mixing dynamical system, central limit theorem leads to the bound $\big|\frac{d \bar{J}^\infty}{ds} - \frac{d \bar{J}}{ds}\big|_{\text{true}}\big| = \mathcal{O}(\frac{1}{\sqrt{T}})$, the total error \cref{eq:sensitivity_error_total} can be expressed as
\begin{equation}
\label{eq:convergence_orders}
   \bigg| \frac{d\bar{J}}{ds} - \frac{d\bar{J}}{ds}\bigg|_{\text{true}}\bigg| \leq \frac{\|e_h\|_2 \|f_s\|_{\infty}}{\sqrt{T}} + \frac{\|f_s\|_\infty \int_0^T \|e_p(t)\| dt}{T} + \mathcal{O}\left(\frac{1}{\sqrt{T}}\right).
\end{equation}
where $\|e_h\|_2$ is bounded because $e_h \in H^1$. With $e_p(t)$ given by \cref{eq:ep_exp}, it can be seen that $\lim\limits_{T\to\infty}\int_0^T\|e_p(t)\| dt < \infty$. Thus, we have $| \frac{d\bar{J}}{ds} - \frac{d\bar{J}}{ds}|_{\text{true}}| \xrightarrow{T\to\infty}0$ and the adjoint LSS sensitivity converges to the true sensitivity. The asymptotic convergence rate of the sensitivity is $\mathcal{O}(\frac{1}{\sqrt{T}})$, which is the same as in \cite{chater_wang_2017} for forward LSS sensitivity derived using a different method.   
\end{proof}
\begin{remark}
\label{rmk:convergence_orders}
  From \cref{eq:convergence_orders}, the error is expected to converge initially at $\mathcal{O}(\frac{1}{T})$ and then as $\mathcal{O}(\frac{1}{\sqrt{T}})$ as $T$ increases.  
\end{remark}

\begin{proposition}
    \label{prop:adjoint_bc}
    If the boundary conditions of the adjoint LSS equation in \cref{eq:adjoint_eq_pre} are non-homogeneous and bounded, the adjoint LSS sensitivity converges to the true sensitivity at the order governed by \cref{eq:convergence_orders}. 
\end{proposition}
\begin{proof}
    Let $\max\{\|\psi(0)\|,\|\psi(T)\|\} < \delta < \infty$. Then, $\|e(0)\| = \|\psi(0) - \psi^\infty (0)\|< \|\psi^\infty (0)\| + \delta < \infty$ and $\|e(T)\| < \infty$. The statement follows by using the modified $e(0)$ and $e(T)$ in \cref{lm:error_bvp} and the derivations that follow until \cref{thm:sensitivity_convergence_rate}. 
\end{proof}

From \cref{lm:L_inverse_bound}, $\|\mathcal{L}^{-1}\|$ is bounded. Thus, it is valid to take a consistent discretization of $\mathcal{L}$, denoted $\mathcal{L}_H$, that is invertible and has a bounded condition number. In fact, it can be proven that certain finite difference schemes are well-conditioned if the continuous differential equation is well-conditioned (c.f. \cite[Theorem 10.75]{ascher_bvp}).   
\begin{corollary}
\label{cor:discrete_truncation_convergence}
    Let $\|\mathcal{L}_H^{-1}\| \leq \kappa$ for mesh sizes of $H < H^*$. Then, for a fixed $T$, the error in the sensitivity due to discretization converges at the order of the local truncation error in the discretization of $\mathcal{L}$.   
\end{corollary}
\begin{proof}
    Let $\psi$ be the solution of the continuous adjoint LSS equation, $\mathcal{L}\psi = g(t)$, and let $\psi_H$ be the discrete solution obtained from the operator $\mathcal{L}_H$. Then, $\|\psi(t)-\psi_H(t)\| = \mathcal{O}(\tau_H)$ \cite{ascher_bvp}, where the truncation error $\tau_H$ is 
    \begin{equation}
        \mathcal{L}_H \psi - g = \tau_H.
    \end{equation}
    The above fact can be seen by noting that, for a well-conditioned system,
\begin{equation}
\label{eq:local_truncation_error}
    \mathcal{L}_H (\psi - \psi_H) = \tau_H \implies \|\psi(t) - \psi_H(t)\| = \|\mathcal{L}_H^{-1}\tau_H\| \leq \kappa \|\tau_H\|, 
\end{equation}
with $\psi(0)-\psi_H(0) = \psi(T)-\psi_H(T)=0$. Assuming that the integration to compute $\frac{d\bar{J}}{ds}$ is discretized with the same order of truncation error $\mathcal{O}(\tau_H)$ as $\mathcal{L}_H$ yields
\begin{equation}
\label{eq:discretization_error_sensitivity}
    \bigg|\frac{d\bar{J}}{ds} - \frac{d\bar{J}}{ds}\bigg|_H\bigg| = \bigg|\frac{1}{T} \int_H (\psi-\psi_H)^T f_s dt + \mathcal{O}(\tau_H)\bigg| \leq \|f_s\|_\infty \kappa \|\tau_H\| + \mathcal{O}(\tau_H) = \mathcal{O}(\tau_H) 
\end{equation}
where $\int_H$ is the discretization of $\int_0^T$. 
\end{proof}

\subsection{Discussion}
In summary, the present paper has shown that the absolute condition number of the LSS operator, $\|\mathcal{L}^{-1}\|$ is bounded as $T\to\infty$ (\cref{lm:L_inverse_bound}). Using \cref{lm:L_inverse_bound}, it was proven that the error $e(t)$ lies in the Sobolev space (\cref{cor:eh_is_in_sobolev}), which led to a proof of existence and uniqueness of the solution to the adjoint LSS and an alternate proof of the convergence of LSS sensitivity to the true sensitivity (\cref{thm:sensitivity_convergence_rate}). It was also shown that the non-homogeneous boundary conditions of the adjoint LSS still lead to the LSS sensitivity converging to the true sensitivity at the rate governed by \cref{eq:convergence_orders} and that the error in the sensitivity due to discretization is of the order of the local truncation error (\cref{cor:discrete_truncation_convergence}). The next section presents test cases to verify the results.

\section{Numerical test cases}
\label{sec:numerical_cases}
Test cases of applications of the method are presented. In particular, we show that the adjoint LSS sensitivity converges at the rate noted in \cref{eq:convergence_orders} irrespective of the adjoint boundary conditions and that the coarse discretized adjoint sensitivity converges at the rate of local truncation error of discretization. 

In the test cases that follow, we take the discretization of time using a constant step size $\Delta t = \frac{T}{N}$ and the temporal grid as $t_i = i\Delta t$.  \cref{eq:lss_adjoint} can then be discretized with a second-order finite difference:
\begin{equation}
\label{eq:second_order_discretization}
    \begin{aligned}
         \frac{\psi_{i+1} - \psi_i}{\Delta t} + f^*_{u_{i+\frac{1}{2}}} \left(\frac{\psi_{i+1} + \psi_i}{2} \right) + j_{u_{i+\frac{1}{2}}} &= r_{i+\frac{1}{2}} \quad i=0,\ldots, N-1 \\
         \frac{r_{i+\frac{1}{2}}-r_{i-\frac{1}{2}}}{\Delta t} - \frac{1}{2} \left(f_{u_{i+\frac{1}{2}}} r_{i+\frac{1}{2}} + f_{u_{i-\frac{1}{2}}} r_{i-\frac{1}{2}}\right)
          &= \frac{1}{\alpha^2}\left( \psi_i^T f_i +  \frac{j_{i+\frac{1}{2}} + j_{i-\frac{1}{2}}}{2}  - \bar{j} \right)f_i \\ & \qquad\qquad \qquad \qquad i=1,\ldots, N-1.  
    \end{aligned}
\end{equation}

with $f_i = \frac{u_{1+\frac{1}{2}} - u_{i-\frac{1}{2}}}{\Delta t}$. The time-averaged integral used for sensitivity analysis in \cref{eq:adjoint_sensitivity} is discretized using the second-order trapezoidal discretization:
\begin{equation}
    \frac{d\bar{J}}{ds} = \frac{1}{N} \sum_{i=0}^{N-1} \left( \left(\frac{\psi_{i+1} + \psi_i}{2} \right)^T f_s|_{i+\frac{1}{2}} + j_s|_{i+\frac{1}{2}}\right).
\end{equation}
A first order finite difference discretization is also used:
\begin{equation}
\label{eq:first_order_discretization}
    \begin{aligned}
         &\frac{\psi_{i+1} - \psi_i}{\Delta t} + f^*_{u_i}\psi_i + j_{u_i} = r_i \quad i=0,\ldots N-1 \\
         &\frac{r_i-r_{i-1}}{\Delta t} - f_{u_i}r_i = \frac{1}{\alpha^2} \left( \psi_i^T f_i + j_i - \bar{j} \right)f_i \quad i=1, \ldots N
    \end{aligned}
\end{equation}

Note that, after eliminating $r_{i}$, one obtains the same LSS discretization matrix as in \cite{wang_lss} for \cref{eq:second_order_discretization} and a similar matrix for \cref{eq:first_order_discretization}. 

\subsection{Lorenz 63}
\begin{figure}[H]
    \centering
    \includegraphics[scale=0.4]{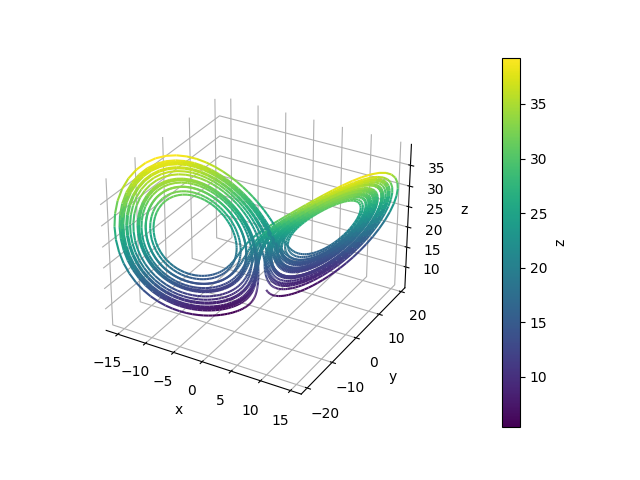}
    \caption{Lorenz attractor}
    \label{fig:Lorenz_attractor}
\end{figure}

Consider the evolution of the state $(x,y,z)$ governed by the Lorenz 63 system which is used to model atmospheric convection \cite{lorentz_63}:
\begin{equation}
    \frac{d}{dt} 
    \begin{bmatrix}
    x \\
    y \\
    z
    \end{bmatrix} =
    \begin{bmatrix}
        \sigma (y-x) \\
        x(\rho - (z-z_0)) - y \\
        xy - \beta(z-z_0)
    \end{bmatrix}
\end{equation}
with $z_0 = 0, \sigma=10, \beta=\frac{8}{3}$ and $\rho=25$. With these parameters, the flow is chaotic with a quasi-hyperbolic strange attractor \cite{sparrow_lorenz}, as shown in \cref{fig:Lorenz_attractor}. 
\begin{figure}[H]
    \centering
    \includegraphics[scale=0.4]{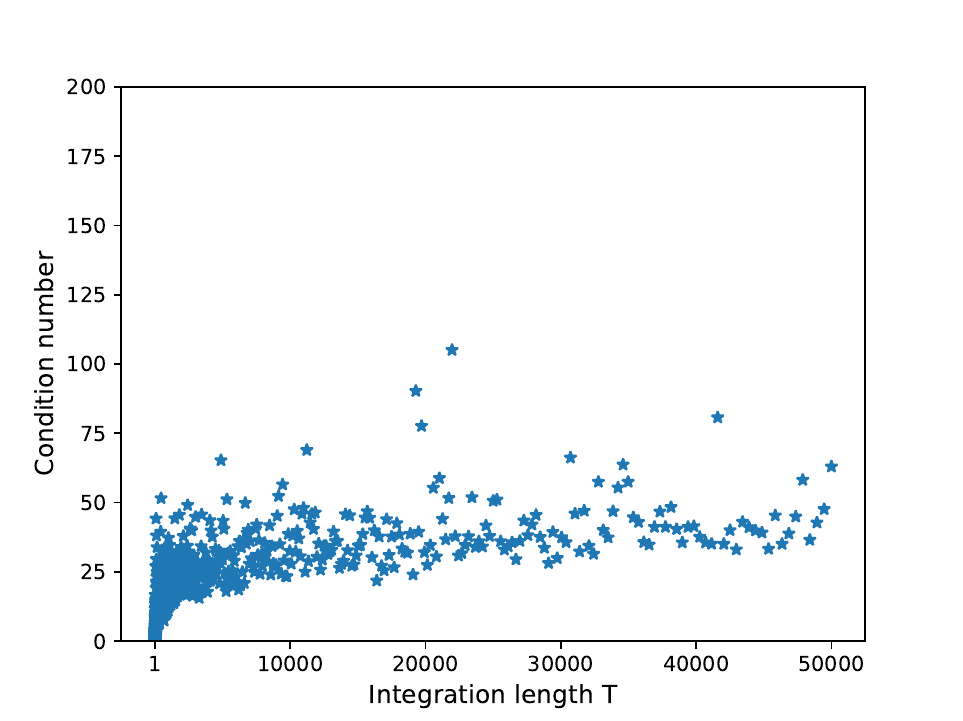}
    \caption{Variation of the LSS condition number, $\|\mathcal{L}^{-1}\|$, with time.}
    \label{fig:condition_number}
\end{figure}
We are interested in the time average of the state $z(t)$:
\begin{equation}
    \bar{J} = \lim_{T\to\infty} \frac{1}{T} \int_0^T z(t) dt,
\end{equation}
and the parameter of interest is $z_0$. As noted in \cite{chater_wang_2017}, changing $z_0$ translates the entire attractor by $z_0$ along the $z$-axis. Thus, the true sensitivity is $\frac{d\bar{J}}{dz_0} = 1$. 

As noted in \cref{rm:L_inverse_l2}, the bound of $\|\mathcal{L}^{-1} \|$ also holds in the 2-norm. In order to verify that the condition number is bounded with $T$ as predicted by \cref{lm:L_inverse_bound}, \cref{fig:condition_number} plots $\|\mathcal{L}^{-1}\| \approx \frac{1}{\lambda_{\min}(\mathcal{L}_H)}$ vs $T$ with $\alpha^2=1$ and for a fixed discretization of $\Delta t = 0.02$, where $\lambda_{\min}(\mathcal{L}_H)$ is the minimum eigenvalue of $\mathcal{L}_H$. As seen from \cref{fig:condition_number}, the condition number remains bounded for large $T$, consistent with \cref{lm:L_inverse_bound}. The boundedness of the condition number was also shown in the test case in \cite{wang_lss}.  

\begin{figure}[H]
    \centering
    \begin{subfigure}[b]{0.49\textwidth}
    \centering
    \includegraphics[width=\textwidth]{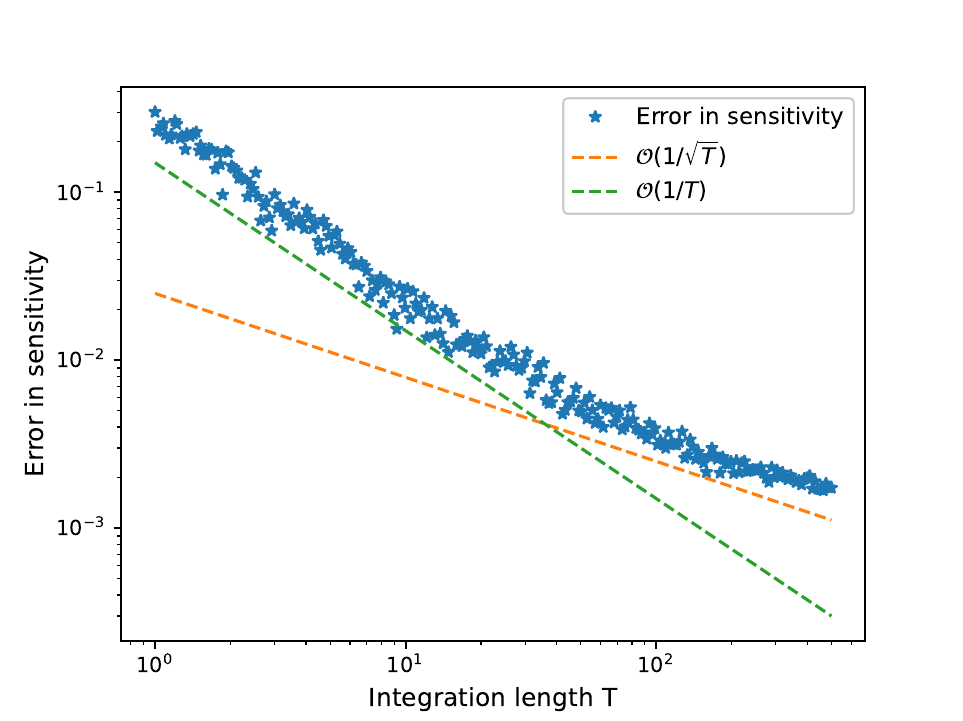}
    \caption{}
    \label{fig:hom_adjoint}      
    \end{subfigure}    
    \hfill
    \begin{subfigure}[b]{0.49\textwidth}
    \centering
    \includegraphics[width=\textwidth]{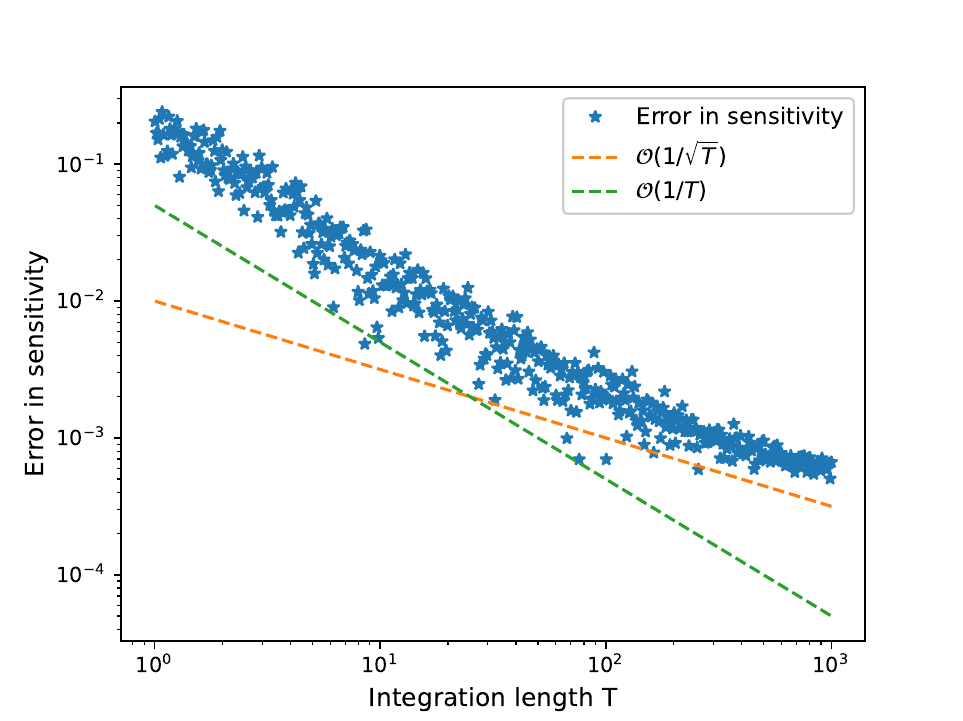}
    \caption{}
    \label{fig:nonhom_adj_bc}      
    \end{subfigure}
    \caption{Convergence of sensitivity with (a) homogeneous and (b) non-homogeneous adjoint boundary conditions.}
    \label{fig:lorentzconvergence}
\end{figure}

To verify the rate of convergence, the adjoint LSS equation is first run by using $\psi(0)=\psi(T)=0$. Therefore, the sensitivity computed is equivalent to that obtained from the forward sensitivity. To verify the effect of integration time $T$, we use the same discrete time domain for both the primal and the adjoint equations, using $\Delta t=0.02$ and $\alpha^2 = 100$. The adjoint LSS algorithm is run and the plot of the error in sensitivity, $\bigg|\frac{d\bar{J}}{ds} - \frac{d\bar{J}}{ds}\bigg|_{\text{true}} \bigg|$, versus integration length $T$ is shown in \cref{fig:hom_adjoint}. The sensitivity at each integration time was calculated by averaging the sensitivities over 20 trajectories obtained with random initial conditions. As seen in \cref{fig:hom_adjoint}, the sensitivity converges at $\mathcal{O}\left(\frac{1}{T}\right)$ initially and later as $\mathcal{O}\left(\frac{1}{\sqrt{T}}\right)$, as noted in \cref{rmk:convergence_orders}.   

Next, to verify that the adjoint boundary conditions do not impact the convergence of the adjoint LSS sensitivity, we set $\psi(0) = \psi(T) = \frac{1}{4}\vect{1}$, where $\vect{1}$ is a vector of ones and $\Delta t = 0.015$. \Cref{fig:nonhom_adj_bc} shows the convergence of the error in the sensitivity with integration time and that it is consistent with \cref{prop:adjoint_bc}; where the error converges at the expected order.

\begin{figure}[H]
    \centering
    \includegraphics[scale=0.4]{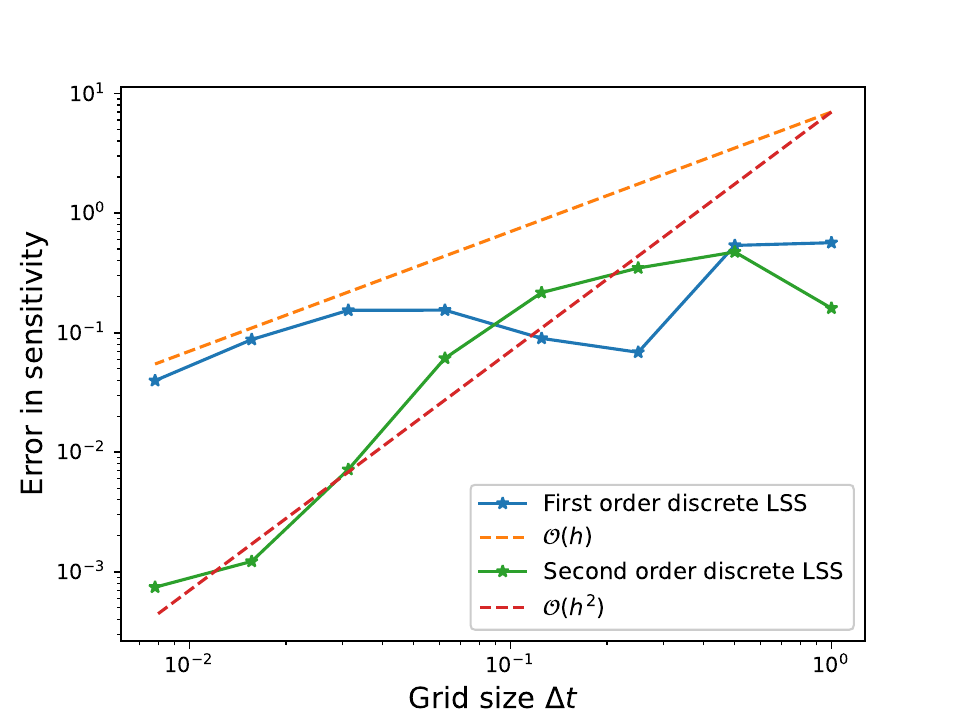}
    \caption{Convergence of sensitivity with grid size using adjoint LSS.}
    \label{fig:discrete_convergence}
\end{figure}

Finally, we verify the impact of the discretization of the adjoint LSS sensitivity and, in particular, the possibility of using a coarse time discretization to reduce the size of the space-time system. Existence and uniqueness of the solution to adjoint LSS, \cref{thm:existence_uniqueness}, guarantees that the discrete adjoint solution obtained from a consistent and stable scheme will converge to the continuous LSS adjoint solution as $h\to 0$. For this run, we use a fixed integration length of $T = 400$ and homogeneous adjoint boundary conditions. \cref{fig:discrete_convergence} shows the convergence of the error in sensitivity with grid size $H$ for both first and second order discretization. The figure shows that the sensitivity error converges at the rate of the local truncation error of the scheme used to discretize the adjoint LSS, \cref{cor:discrete_truncation_convergence}.

\subsection{Coupled oscillator}
The adjoint LSS was applied on the model A coupled oscillator of Kuznetsov and Pikovsky \cite{kuznetsov}, which was numerically shown to be hyperbolic. 
\begin{figure}[H]
    \centering
    \includegraphics[scale=0.4]{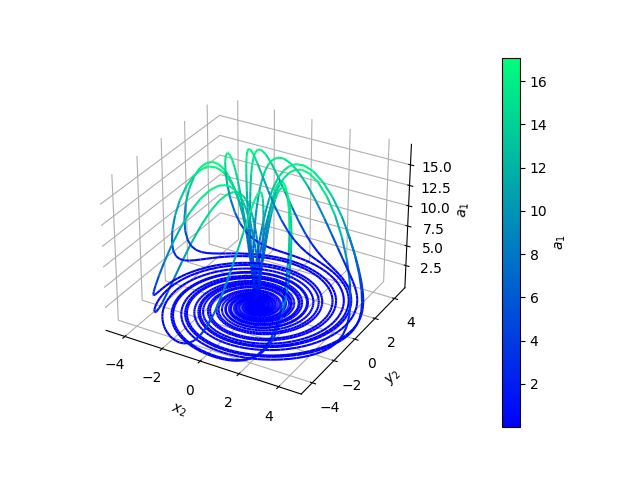}
    \caption{Attractor of the coupled oscillator}
    \label{fig:oscillator_attractor}
\end{figure}
The oscillator has the form  
\begin{align*}
     & \dot{x_1} = \omega_0 y_1 + \left(1-a_2+\frac{1}{2}a_1-\frac{1}{50}a_1^2\right) x_1 + \epsilon (x_2-s) y_2 \\
     & \dot{y_1} = -\omega_0 x_1 + \left(1-a_2+\frac{1}{2}a_1-\frac{1}{50}a_1^2\right) y_1 \\
     & \dot{x_2} = \omega_0 y_2 + \left(a_1-1\right) (x_2-s) + \epsilon x_1 \\
     & \dot{y_2} = -\omega_0 (x_2-s) + \left(a_1-1\right) y_2 
\end{align*}
with $\omega_0 = 2\pi$, $\epsilon=0.3$, $s=0$, $a_1 = x_1^2 + y_1^2$ and $a_2 = (x_2-s)^2 + y_2^2$. \Cref{fig:oscillator_attractor} shows the attractor of the oscillator.

We are interested in evaluating the sensitivity of the following functional with respect to the parameter $s$:
\begin{equation}
      \bar{J} = \lim_{T\to\infty} \frac{1}{T} \int_0^T x_2(t) dt.
\end{equation}
\begin{figure}[H]
    \centering
    \includegraphics[scale=0.4]{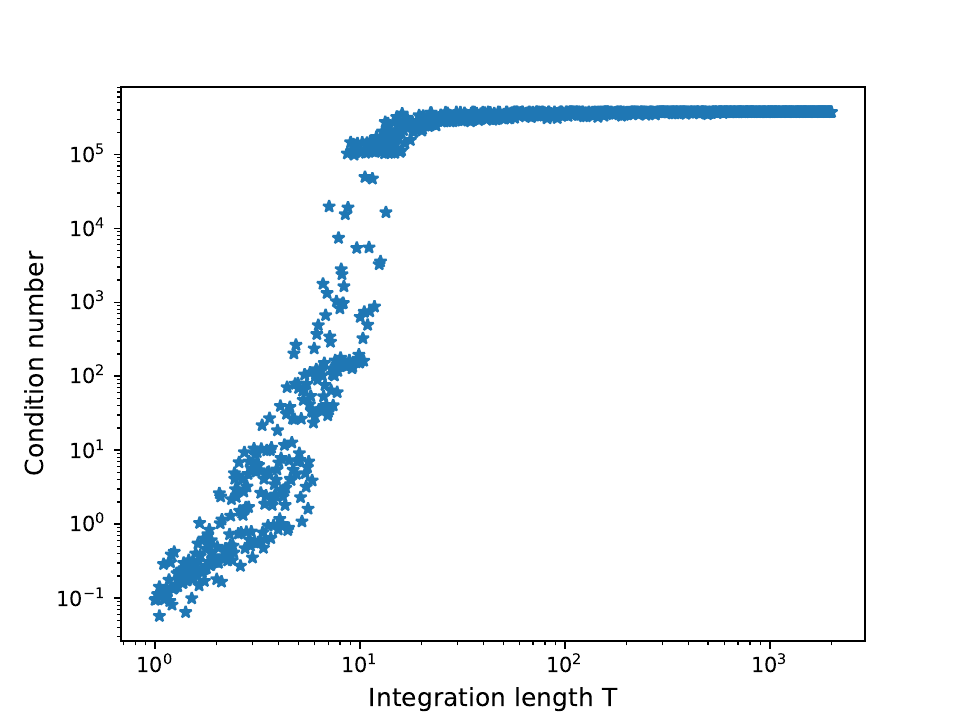}
    \caption{Variation of the absolute LSS condition number, $\|\mathcal{L}^{-1}\|$, with integration time.}
    \label{fig:condition_number_oscillator}
\end{figure}

Similar to the previous test case, the exact value of the sensitivity is known to be 1. Adjont LSS was applied with a time step of $\Delta t=0.01$ and $\alpha^2=100$ for various integration lengths. \Cref{fig:condition_number_oscillator} shows the variation of the condition number with the integration time $T$. The absolute condition number, $\|\mathcal{L}_H^{-1}\|$, though quite large, is bounded after a certain time and does not scale with integration length $T$ in accordance with \cref{lm:L_inverse_bound}.  

\begin{figure}[H]
    \centering
    \begin{subfigure}[b]{0.49\textwidth}
    \centering
    \includegraphics[width=\textwidth]{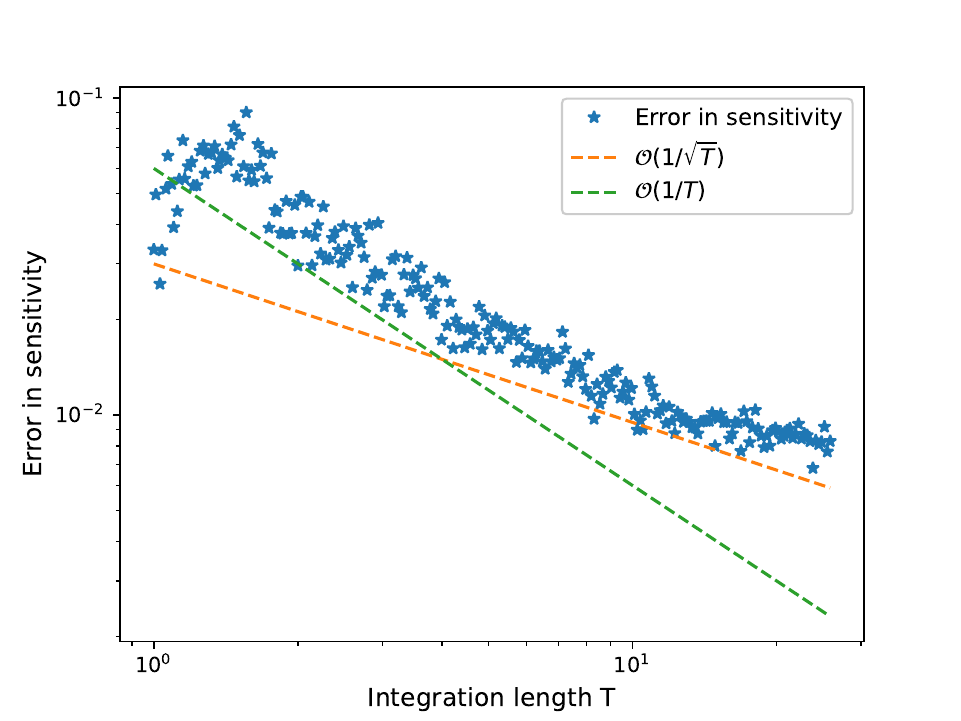}
    \caption{}
    \label{fig:oscillator_hom_adjoint}      
    \end{subfigure}    
    \hfill
    \begin{subfigure}[b]{0.49\textwidth}
    \centering
    \includegraphics[width=\textwidth]{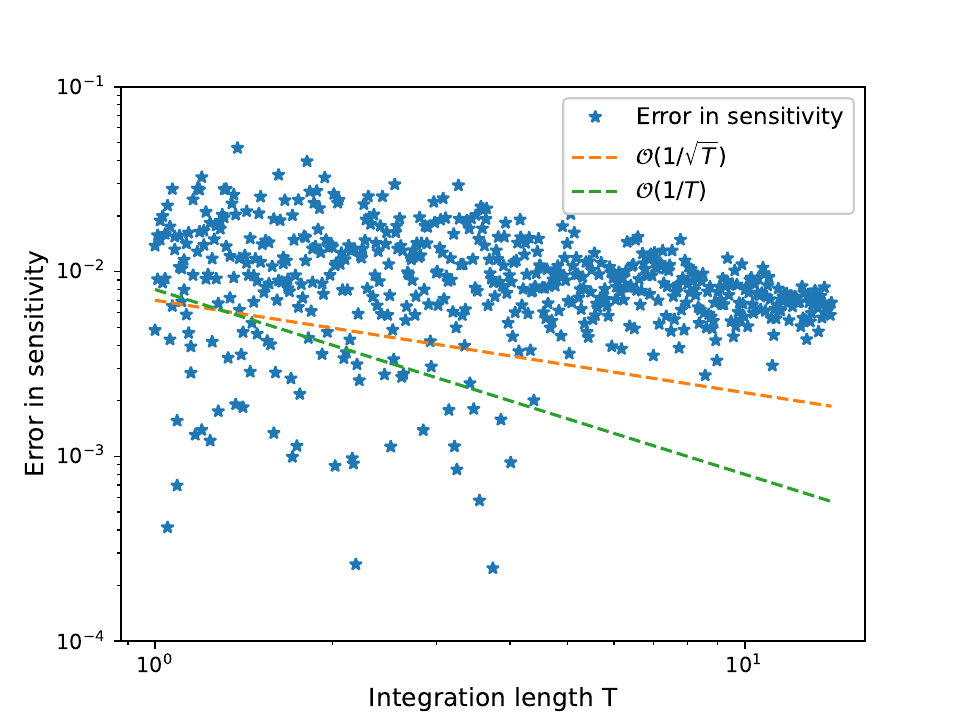}
    \caption{}
    \label{fig:oscillator_nonhom_adjoint}      
    \end{subfigure}
    \caption{Convergence of sensitivity with (a) homogeneous and (b) non-homogeneous adjoint boundary conditions.}
    \label{fig:oscillatorconvergence}
\end{figure}

Adjoint LSS was run using homogeneous adjoint boundary conditions with $\Delta t=0.01$ and $\alpha^2=100$. The sensitivity value for each integration time was obtained by averaging the sensitivities from 20 trajectories starting from random initial conditions with a spin-up time of $100$ units to get the solution on the attractor. \Cref{fig:oscillator_hom_adjoint} plots the error in the sensitivities with integration time. The error converges initially as $\mathcal{O}\left(\frac{1}{T}\right)$ and then as $\mathcal{O}\left(\frac{1}{\sqrt{T}}\right)$ as predicted by \cref{rmk:convergence_orders}. Adjoint LSS was also run using non-homogeneous boundary condition of $\psi = \frac{1}{5} \mathbf{1}$ with $\alpha^2=1$. \Cref{fig:oscillator_nonhom_adjoint} shows the convergence of the error in sensitivity with the integration time. We note that the variation in sensitivities is higher in this case, with a visible $\mathcal{O}\left(\frac{1}{\sqrt{T}}\right)$ convergence rate. 

\begin{figure}[H]
    \centering
    \includegraphics[scale=0.4]{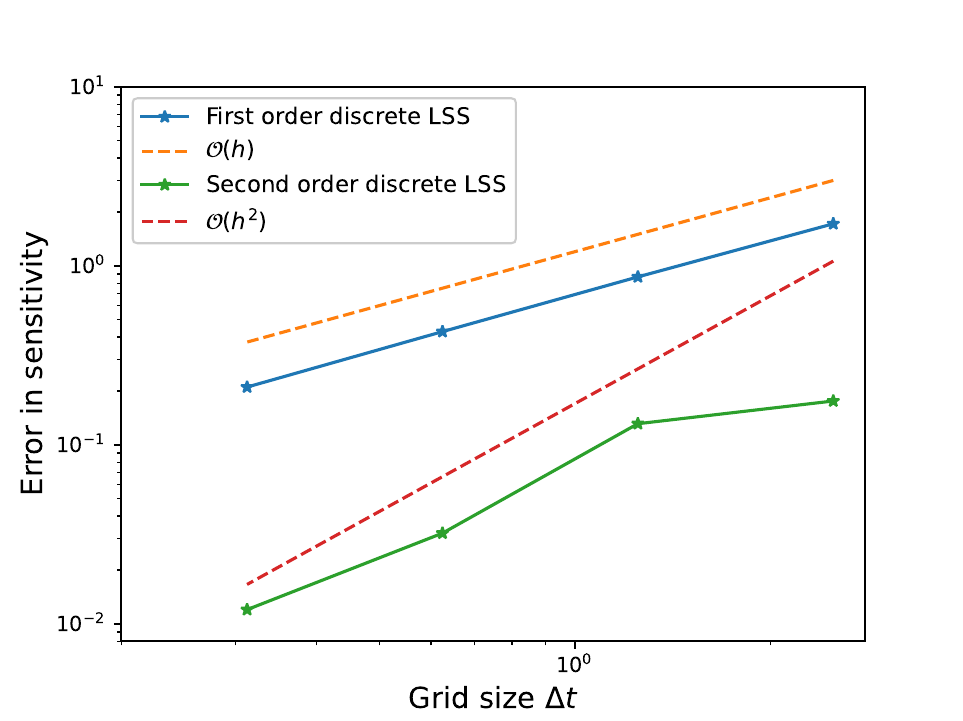}
    \caption{Convergence of sensitivity with grid size using adjoint LSS.}
    \label{fig:oscillator_dt_convergence}
\end{figure}

To evaluate the impact of error due to discretization, first- and second-order accurate adjoint LSS was run with a coarse temporal discretization, using $\alpha^2 = 1$ with averaging over $10$ random trajectories. The trajectories of the solution were computed using $\Delta t=0.01$ before interpolating the solution to a coarse grid. \Cref{fig:oscillator_dt_convergence} shows that the error in sensitivity converges at the expected order of the scheme's truncation error. 
\section{Conclusion}
An analysis of the adjoint of the least squares shadowing (LSS) approach was conducted. It was found that the condition number of the LSS equation is bounded for large integration times in the case of uniformly hyperbolic systems. It was also found that the error in the adjoint solution lies in the Sobolev space $H^1$. Hence, the norms of both the adjoint error and the time derivative of error are bounded. The existence and uniqueness of the solution to the adjoint LSS were shown, along with an alternate proof of convergence of the LSS to the true sensitivity. It was also shown that modifying the adjoint boundary conditions does not affect the rate of convergence of the method. The bounds on the conditioning of the LSS derived in the current work indicate the trend between the conditioning and the time dilation factor, which is consistent with the observation in numerical cases in previous LSS literature. The existence and uniqueness of the solution to the continuous adjoint LSS ensures that any stable and consistent temporal discretization of the adjoint LSS would yield a solution that converges to the continuous adjoint solution as the time step is refined. We showed that the condition number of the LSS operator is bounded and numerically verified the result. It was also shown that the error due to the temporal discretization of the LSS, scales at the order of the local truncation error of the discretization scheme. 

Several issues need to be addressed to be able to use adjoint LSS for large-scale problems, particularly for problems involving several positive Lyapunov exponents for which the alternate approaches of computing the unstable subspaces, such as NILSS  and NILSAS \cite{nilss, nilsas}, can be expensive. In the present work, the size of the LSS matrix was reduced by solving the adjoint LSS using a coarse discretization of time. This approach can be improved by determining an optimal lower-dimensional subspace to solve the LSS problem. The errors introduced due to modifications in continuous adjoint LSS were shown to be bounded by local operations of the approximation. Hence, the error due to such approximations is expected to be low, though additional investigation is necessary to determine the magnitude of approximation errors a-priori. In \cite{chandramoorthy_wang_nonphysical} it was found that shadowing-based approaches can yield sensitivities with bias if the shadowing trajectory is not representative of a true physical trajectory for cases that are not uniformly hyperbolic. Whether this problem is confined to lower-dimensional test cases or also applies to problems in aerodynamics with large number of degrees of freedom needs further investigation.

\bibliographystyle{siamplain}

\begin{thebibliography}{10}

\bibitem{optimal_control}
{\sc T.~Akman, B.~Hosseini, J.~Diepolder, B.~Grüter, R.~Afonso, M.~Gerdts, and F.~Holzapfel}, {\em Efficient sensitivity calculation for robust optimal control}, 09 2020, \url{https://doi.org/10.25967/490248}.

\bibitem{anosov_shadowing}
{\sc D.~V. Anosov}, {\em Geodesic flows on closed {R}iemannian manifolds of negative curvature}, Trudy Matematicheskogo Instituta Imeni VA Steklova, 90 (1967), pp.~3--210.

\bibitem{ascher_bvp}
{\sc U.~M. Ascher, R.~M. Mattheij, and R.~D. Russell}, {\em Numerical solution of boundary value problems for ordinary differential equations}, Society for Industrial and Applied Mathematics, 1995.

\bibitem{ashley_hicken_fourier}
{\sc A.~S. Ashley and J.~E. Hicken}, {\em Sensitivity Analysis of Chaotic Problems using a Fourier Approximation of the Least-Squares Adjoint}, \url{https://doi.org/10.2514/6.2016-4409}, \url{https://arxiv.org/abs/https://arc.aiaa.org/doi/pdf/10.2514/6.2016-4409}.

\bibitem{becker_rannacher_2001}
{\sc R.~Becker and R.~Rannacher}, {\em An optimal control approach to a posteriori error estimation in finite element methods}, Acta Numerica, 10 (2001), p.~1–102, \url{https://doi.org/10.1017/S0962492901000010}.

\bibitem{blonigan_multigrid}
{\sc P.~Blonigan and Q.~Wang}, {\em Multigrid-in-time for sensitivity analysis of chaotic dynamical systems}, Numerical Linear Algebra with Applications, 24 (2014), p.~e1946, \url{https://doi.org/https://doi.org/10.1002/nla.1946}, \url{https://arxiv.org/abs/https://onlinelibrary.wiley.com/doi/pdf/10.1002/nla.1946}.

\bibitem{blonigan_adjoint_nilss}
{\sc P.~J. Blonigan}, {\em Adjoint sensitivity analysis of chaotic dynamical systems with non-intrusive least squares shadowing}, Journal of Computational Physics, 348 (2017), pp.~803--826, \url{https://doi.org/https://doi.org/10.1016/j.jcp.2017.08.002}.

\bibitem{blonigan_fokker_planck}
{\sc P.~J. Blonigan and Q.~Wang}, {\em Probability density adjoint for sensitivity analysis of the mean of chaos}, Journal of Computational Physics, 270 (2014), pp.~660--686, \url{https://doi.org/https://doi.org/10.1016/j.jcp.2014.04.027}.

\bibitem{blonigan_mss}
{\sc P.~J. Blonigan and Q.~Wang}, {\em Multiple shooting shadowing for sensitivity analysis of chaotic dynamical systems}, Journal of Computational Physics, 354 (2018), pp.~447--475, \url{https://doi.org/https://doi.org/10.1016/j.jcp.2017.10.032}.

\bibitem{blonigan_airfoil}
{\sc P.~J. Blonigan, Q.~Wang, E.~J. Nielsen, and B.~Diskin}, {\em Least-squares shadowing sensitivity analysis of chaotic flow around a two-dimensional airfoil}, AIAA Journal, 56 (2018), pp.~658--672, \url{https://doi.org/10.2514/1.J055389}, \url{https://arxiv.org/abs/https://doi.org/10.2514/1.J055389}.

\bibitem{bowen_shadowing}
{\sc R.~Bowen}, {\em $\omega$-limit sets for {Axiom A} diffeomorphisms}, Journal of Differential Equations, 18 (1975), pp.~333--339, \url{https://doi.org/https://doi.org/10.1016/0022-0396(75)90065-0}.

\bibitem{brenner_scott}
{\sc S.~C. Brenner and L.~Scott}, {\em The mathematical theory of finite element methods}, Springer, 2008.

\bibitem{wang_ensemble_adjoint}
{\sc N.~Chandramoorthy, P.~Fernandez, C.~Talnikar, and Q.~Wang}, {\em An Analysis of the Ensemble Adjoint Approach to Sensitivity Analysis in Chaotic Systems}, \url{https://doi.org/10.2514/6.2017-3799}, \url{https://arxiv.org/abs/https://arc.aiaa.org/doi/pdf/10.2514/6.2017-3799}.

\bibitem{chandramoorthy_wang_nonphysical}
{\sc N.~Chandramoorthy and Q.~Wang}, {\em On the probability of finding nonphysical solutions through shadowing}, Journal of Computational Physics, 440 (2021), p.~110389, \url{https://doi.org/https://doi.org/10.1016/j.jcp.2021.110389}.

\bibitem{chandramoorthy_efficient_linear_response}
{\sc N.~Chandramoorthy and Q.~Wang}, {\em Efficient computation of linear response of chaotic attractors with one-dimensional unstable manifolds}, SIAM Journal on Applied Dynamical Systems, 21 (2022), pp.~735--781, \url{https://doi.org/10.1137/21M1405599}, \url{https://arxiv.org/abs/https://doi.org/10.1137/21M1405599}.

\bibitem{chater_wang_2017}
{\sc M.~Chater, A.~Ni, P.~J. Blonigan, and Q.~Wang}, {\em Least squares shadowing method for sensitivity analysis of differential equations}, SIAM Journal on Numerical Analysis, 55 (2017), pp.~3030--3046, \url{https://doi.org/10.1137/15M1039067}, \url{https://arxiv.org/abs/https://doi.org/10.1137/15M1039067}.

\bibitem{corrigan_et_al_2019}
{\sc A.~Corrigan, A.~D. Kercher, and D.~A. Kessler}, {\em A moving discontinuous {Galerkin} finite element method for flows with interfaces}, International Journal for Numerical Methods in Fluids, 89 (2019), pp.~362--406.

\bibitem{dolejsigoal1}
{\sc V.~Dolejší, O.~Bartoš, and F.~Roskovec}, {\em Goal-oriented mesh adaptation method for nonlinear problems including algebraic errors}, Computers \& Mathematics with Applications, 93 (2021), pp.~178--198.

\bibitem{eckmann_ergodic}
{\sc J.-P. Eckmann and D.~Ruelle}, {\em Ergodic theory of chaos and strange attractors}, Reviews of modern physics, 57 (1985), p.~617.

\bibitem{lea_ensemble_adjoint}
{\sc G.~L. Eyink, T.~W.~N. Haine, and D.~J. Lea}, {\em Ruelle's linear response formula, ensemble adjoint schemes and l{\'e}vy flights}, Nonlinearity, 17 (2004), pp.~1867--1889.

\bibitem{fidkowski_darmofal_2011}
{\sc K.~J. Fidkowski and D.~L. Darmofal}, {\em Review of output-based error estimation and mesh adaptation in computational fluid dynamics}, AIAA Journal, 49 (2011), pp.~673--694, \url{https://doi.org/10.2514/1.J050073}, \url{https://arxiv.org/abs/https://doi.org/10.2514/1.J050073}.

\bibitem{giles_suli}
{\sc M.~B. Giles and E.~Süli}, {\em Adjoint methods for {PDE}s: a posteriori error analysis and postprocessing by duality}, Acta Numerica, 11 (2002), p.~145–236, \url{https://doi.org/10.1017/S096249290200003X}.

\bibitem{ginelli_2007}
{\sc F.~Ginelli, P.~Poggi, A.~Turchi, H.~Chat\'e, R.~Livi, and A.~Politi}, {\em Characterizing dynamics with covariant lyapunov vectors}, Phys. Rev. Lett., 99 (2007), p.~130601, \url{https://doi.org/10.1103/PhysRevLett.99.130601}.

\bibitem{numerical_shadowing}
{\sc S.~M. Hammel, J.~A. Yorke, and C.~Grebogi}, {\em {Numerical orbits of chaotic processes represent true orbits}}, Bulletin (New Series) of the American Mathematical Society, 19 (1988), pp.~465--469.

\bibitem{hicken_chaos}
{\sc J.~Hicken and V.~Ramakrishnan}, {\em A method to regularize optimization problems governed by chaotic dynamical systems}, Chaos, Solitons \& Fractals, 188 (2024), p.~115491, \url{https://doi.org/https://doi.org/10.1016/j.chaos.2024.115491}.

\bibitem{jameson_adjoint_optimization}
{\sc A.~Jameson}, {\em Aerodynamic design via control theory}, Journal of scientific computing, 3 (1988), pp.~233--260.

\bibitem{kantarakias_papadakis}
{\sc K.~D. Kantarakias and G.~Papadakis}, {\em Sensitivity analysis of chaotic systems using a frequency-domain shadowing approach}, Journal of Computational Physics, 474 (2023), p.~111757, \url{https://doi.org/https://doi.org/10.1016/j.jcp.2022.111757}.

\bibitem{vermier_lss}
{\sc H.~R. Karbasian and B.~C. Vermeire}, {\em {Sensitivity analysis of chaotic dynamical systems using a physics-constrained data-driven approach}}, Physics of Fluids, 34 (2022), p.~014101, \url{https://doi.org/10.1063/5.0076074}, \url{https://arxiv.org/abs/https://pubs.aip.org/aip/pof/article-pdf/doi/10.1063/5.0076074/16618988/014101\_1\_online.pdf}.

\bibitem{katok_dynamical_system}
{\sc A.~Katok and B.~Hasselblatt}, {\em Introduction to the Modern Theory of Dynamical Systems}, Encyclopedia of Mathematics and its Applications, Cambridge University Press, 1995.

\bibitem{adjoint_clv}
{\sc P.~V. Kuptsov and U.~Parlitz}, {\em Theory and computation of covariant lyapunov vectors}, Journal of Nonlinear Science, 22 (2011), pp.~727 -- 762.

\bibitem{kuznetsov}
{\sc S.~P. Kuznetsov and A.~Pikovsky}, {\em Autonomous coupled oscillators with hyperbolic strange attractors}, Physica D: Nonlinear Phenomena, 232 (2007), pp.~87--102, \url{https://doi.org/https://doi.org/10.1016/j.physd.2007.05.008}.

\bibitem{lasagna_upo}
{\sc D.~Lasagna}, {\em Sensitivity analysis of chaotic systems using unstable periodic orbits}, SIAM Journal on Applied Dynamical Systems, 17 (2018), pp.~547--580, \url{https://doi.org/10.1137/17M114354X}, \url{https://arxiv.org/abs/https://doi.org/10.1137/17M114354X}.

\bibitem{lasagna_periodic_shadowing}
{\sc D.~Lasagna, A.~Sharma, and J.~Meyers}, {\em Periodic shadowing sensitivity analysis of chaotic systems}, Journal of Computational Physics, 391 (2019), pp.~119--141, \url{https://doi.org/https://doi.org/10.1016/j.jcp.2019.04.021}.

\bibitem{lea_e_al_2000}
{\sc D.~J. Lea, M.~R. Allen, and T.~W.~N. Haine}, {\em Sensitivity analysis of the climate of a chaotic system}, Tellus A: Dynamic Meteorology and Oceanography,  (2000), \url{https://doi.org/10.3402/tellusa.v52i5.12283}.

\bibitem{lorentz_63}
{\sc E.~N. Lorenz}, {\em Deterministic nonperiodic flow}, Journal of Atmospheric Sciences, 20 (1963), pp.~130 -- 141, \url{https://doi.org/10.1175/1520-0469(1963)020<0130:DNF>2.0.CO;2}.

\bibitem{inriagoaloriented}
{\sc A.~Loseille, A.~Dervieux, and F.~Alauzet}, {\em {Fully anisotropic goal-oriented mesh adaptation for 3D steady Euler equations}}, Journal of Computational Physics, 229 (2010), pp.~2866--2897.

\bibitem{ni_adjoint_arxiv}
{\sc A.~Ni}, {\em Adjoint shadowing directions in hyperbolic systems for sensitivity analysis}, arXiv preprint arXiv:1807.05568,  (2018).

\bibitem{nilsas}
{\sc A.~Ni and C.~Talnikar}, {\em Adjoint sensitivity analysis on chaotic dynamical systems by non-intrusive least squares adjoint shadowing ({NILSAS})}, Journal of Computational Physics, 395 (2019), pp.~690--709, \url{https://doi.org/https://doi.org/10.1016/j.jcp.2019.06.035}.

\bibitem{ni2023recursive}
{\sc A.~Ni and Y.~Tong}, {\em Recursive divergence formulas for perturbing unstable transfer operators and physical measures}, Journal of Statistical Physics, 190 (2023), p.~126.

\bibitem{nilss}
{\sc A.~Ni and Q.~Wang}, {\em Sensitivity analysis on chaotic dynamical systems by non-intrusive least squares shadowing ({NILSS})}, Journal of Computational Physics, 347 (2017), pp.~56--77, \url{https://doi.org/https://doi.org/10.1016/j.jcp.2017.06.033}.

\bibitem{nocedal_wright}
{\sc J.~Nocedal and S.~J. Wright}, {\em Numerical Optimization}, Springer, New York, NY, USA, 2e~ed., 2006.

\bibitem{oseledets}
{\sc V.~I. Oseledec}, {\em A multiplicative ergodic theorem, lyapunov characteristic numbers for dynamical systems}, Transactions of the Moscow Mathematical Society, 19 (1968), pp.~197--231.

\bibitem{pierce_giles}
{\sc N.~A. Pierce and M.~B. Giles}, {\em Adjoint recovery of superconvergent functionals from {PDE} approximations}, SIAM Review, 42 (2000), pp.~247--264, \url{https://doi.org/10.1137/S0036144598349423}, \url{https://arxiv.org/abs/https://doi.org/10.1137/S0036144598349423}.

\bibitem{pilyugin_shadowing}
{\sc S.~Y. Pilyugin}, {\em Shadowing in dynamical systems}, 1999.

\bibitem{ruelle_ergodic}
{\sc D.~Ruelle}, {\em Ergodic theory of differentiable dynamical systems}, Publications Math{\'e}matiques de l'Institut des Hautes {\'E}tudes Scientifiques, 50 (1979), pp.~27--58.

\bibitem{ruelle_diff_of_srb}
{\sc D.~Ruelle}, {\em Differentiation of {SRB} states}, Communications in Mathematical Physics, 187 (1997), pp.~227--241.

\bibitem{fidkowski_lss}
{\sc Y.~S. Shimizu and K.~Fidkowski}, {\em Output-Based Error Estimation for Chaotic Flows Using Reduced-Order Modeling}, \url{https://doi.org/10.2514/6.2018-0826}, \url{https://arxiv.org/abs/https://arc.aiaa.org/doi/pdf/10.2514/6.2018-0826}.

\bibitem{sliwiak2023}
{\sc A.~A. {\'S}liwiak and Q.~Wang}, {\em Approximating the linear response of physical chaos}, Nonlinear Dynamics, 111 (2023), pp.~1835--1869.

\bibitem{sparrow_lorenz}
{\sc C.~Sparrow}, {\em The Lorenz equations: Bifurcations, Chaos, and Strange Attractors}, vol.~269, Springer New York, 1982.

\bibitem{thakur_nadarajah_2024}
{\sc P.~Thakur and S.~Nadarajah}, {\em Adjoint-based goal-oriented implicit shock tracking using full space mesh optimization}, Journal of Computational Physics, 523 (2025), p.~113633, \url{https://doi.org/https://doi.org/10.1016/j.jcp.2024.113633}.

\bibitem{thuburn_fokker_planck}
{\sc J.~Thuburn}, {\em Climate sensitivities via a fokker–planck adjoint approach}, Quarterly Journal of the Royal Meteorological Society, 131 (2005), pp.~73--92, \url{https://doi.org/https://doi.org/10.1256/qj.04.46}, \url{https://arxiv.org/abs/https://rmets.onlinelibrary.wiley.com/doi/pdf/10.1256/qj.04.46}.

\bibitem{wang_2013}
{\sc Q.~Wang}, {\em Forward and adjoint sensitivity computation of chaotic dynamical systems}, Journal of Computational Physics, 235 (2013), pp.~1--13, \url{https://doi.org/https://doi.org/10.1016/j.jcp.2012.09.007}.

\bibitem{wang_gomez_lss}
{\sc Q.~Wang, S.~A. Gomez, P.~J. Blonigan, A.~L. Gregory, and E.~Y. Qian}, {\em {Towards scalable parallel-in-time turbulent flow simulations}}, Physics of Fluids, 25 (2013), p.~110818, \url{https://doi.org/10.1063/1.4819390}, \url{https://arxiv.org/abs/https://pubs.aip.org/aip/pof/article-pdf/doi/10.1063/1.4819390/13380424/110818\_1\_online.pdf}.

\bibitem{wang_lss}
{\sc Q.~Wang, R.~Hu, and P.~Blonigan}, {\em Least squares shadowing sensitivity analysis of chaotic limit cycle oscillations}, Journal of Computational Physics, 267 (2014), pp.~210--224, \url{https://doi.org/https://doi.org/10.1016/j.jcp.2014.03.002}.

\bibitem{zahr1}
{\sc M.~Zahr and P.-O. Persson}, {\em An optimization-based approach for high-order accurate discretization of conservation laws with discontinuous solutions}, Journal of Computational Physics, 365 (2018), pp.~105--134.

\end{thebibliography}

\end{document}